\documentclass[11pt,a4paper]{article}
\usepackage{amsmath}
\usepackage{amssymb}
\usepackage{amsfonts}
\usepackage{amsthm}
\usepackage{relsize}
\usepackage{setspace}
\usepackage[margin=2.5cm]{geometry}

\usepackage{url}
\usepackage{xspace}
\usepackage{mathrsfs}
\usepackage{tocloft}
\usepackage{ulem}
\usepackage{titling}

\usepackage[shortlabels]{enumitem}
\setlist{topsep=0pt}

\usepackage{hyperref}
\usepackage{authblk}
\usepackage[usenames, dvipsnames]{color}

\usepackage{xr}  
\newtheorem{thm}{Theorem}[section]
\newtheorem{lem}[thm]{Lemma}

\newtheorem{prop}[thm]{Proposition}
\newtheorem{cor}[thm]{Corollary}

\theoremstyle{definition}
\newtheorem{defn}[thm]{Definition}

\newtheorem{rem}[thm]{Remark}

\newcommand{\defbold}{\textbf}

\newcommand{\inv}{^{-1}}

\newcommand{\QZ}{\mathrm{QZ}}
\newcommand{\QC}{\mathrm{QC}}

\newcommand{\CC}{\mathrm{C}}
\newcommand{\M}{\mathrm{M}}
\newcommand{\N}{\mathrm{N}}

\newcommand{\Z}{\mathrm{Z}}

\newcommand{\tdlc}{t.d.l.c.\@\xspace}
\newcommand{\hji}{h.j.i.\@\xspace}

\newcommand{\con}{\mathrm{con}}

\newcommand{\ldlat}{\mathcal{LD}}
\newcommand{\LDlat}{\mathrm{LD}}

\newcommand{\lcent}{\mathcal{LC}}
\newcommand{\LCent}{\mathrm{LC}}

\newcommand{\propLD}{(LD)\@\xspace}

\newcommand{\lnorm}{\mathcal{LN}}

\newcommand{\Hom}{\mathrm{Hom}}

\newcommand{\Aut}{\mathrm{Aut}}

\newcommand{\Sym}{\mathrm{Sym}}

\newcommand{\Res}{\mathrm{Res}}

\newcommand{\bN}{\mathbb{N}}
\newcommand{\bP}{\mathbb{P}}
\newcommand{\bQ}{\mathbb{Q}}

\newcommand{\bT}{\mathbb{T}}
\newcommand{\bZ}{\mathbb{Z}}

\newcommand{\mc}{\mathcal}
\newcommand{\ms}{\mathscr}

\newcommand{\triv}{\{1\}}

\newcommand{\cgrp}[1]{\overline{\langle #1 \rangle}}

\newcommand{\grp}[1]{\langle #1 \rangle}
\newcommand{\ol}[1]{\overline{#1}}

\begin{document}

\title{Totally disconnected locally compact groups with just infinite locally normal subgroups}

\preauthor{\large}
\DeclareRobustCommand{\authoring}{
\renewcommand{\thefootnote}{\arabic{footnote}}
\begin{center}Colin D. Reid\textsuperscript{1}\footnotetext[1]{Research supported by ARC grant FL170100032.}
\\ \bigskip
The University of Newcastle, School of Mathematical and Physical Sciences, Callaghan, NSW 2308, Australia.\\
\href{mailto:colin@reidit.net}{colin@reidit.net}
\end{center}
}
\author{\authoring}
\postauthor{\par}

\maketitle

\begin{abstract}
We obtain a characterization of totally disconnected, locally compact groups $G$ with the following property: given a locally normal subgroup $K$ of $G$, then there is an open subgroup of $K$ that is a direct factor of an open subgroup of $G$.  This property is motivated by J. Wilson's structure theory of just infinite groups, and indeed, when $G$ has trivial quasi-centre, the condition turns out to be equivalent to the condition that $G$ is locally isomorphic to a finite direct product of just infinite profinite groups.  In the latter situation we obtain some global structural features of $G$, building on an earlier result of Barnea--Ershov--Weigel and also using tools developed by P.-E. Caprace, G. Willis and the author for studying local structure in totally disconnected locally compact groups.
\end{abstract}

\tableofcontents

\addtocontents{toc}{\protect\setcounter{tocdepth}{1}}

\section{Introduction}

\subsection{Background}

The background and motivation for this article comes primarily from three sources.  First, we recall some of the structure theory of just infinite groups; see \cite{Gri} and \cite{WilNH} for a detailed account.

\begin{defn}A topological group $G$ is \defbold{just infinite} if $G$ is infinite and every nontrivial closed normal subgroup of $G$ has finite index.

A \defbold{hereditarily just infinite} (\hji) group is one for which every finite index open subgroup is just infinite.

A \defbold{branch group} is a compact or discrete group $G$ acting faithfully by automorphisms on a locally finite rooted tree, with root vertex $\varepsilon$, such that for each sphere of vertices $S_n(\varepsilon)$ around the root, $G$ acts transitively on $S_n(\varepsilon)$ and a finite index subgroup of $G$ splits as a direct product $\prod_{v \in S_n(\varepsilon)}R_v$, where $R_v$ fixes all vertices not descended from $v$.
\end{defn}

\begin{thm}[{\cite[Theorem~3]{WilJI}, \cite[Corollary~4.4]{WilNH}\footnote{The cited results do not claim the result for a just infinite virtually abelian profinite group $G$, but this case can be verified by a similar argument to \cite[Theorem~3A]{WilJI}.}}]\label{intro:wilcomp}
Let $G$ be a just infinite profinite or discrete group.  Then every closed subnormal subgroup of $G$ is a direct factor of a finite index subgroup.
\end{thm}

\begin{thm}[{See \cite[Theorem~3]{Gri}}]\label{intro:wiltri}Let $G$ be a just infinite profinite or discrete group that is not virtually abelian.  Then exactly one of the following holds:
\begin{enumerate}[(i)]
\item There is a (closed) normal subgroup of $G$ of finite index of the form $L_1 \times \dots \times L_n$, where the $L_i$ form a single conjugacy class of hereditarily just infinite subgroups of $G$;
\item $G$ is a branch group.
\end{enumerate}
\end{thm}

In the present article we will consider a property inspired by J. Wilson's Theorem~\ref{intro:wilcomp}, converted into a form expressible in the local structure of totally disconnected, locally compact (\tdlc) groups.

\begin{defn}
Let $G$ be a \tdlc group.  Say that $G$ has \defbold{property~\propLD} if for every closed subgroup $K$ of $G$ with open normalizer, then there is an open subgroup $L$ of $K$ that is a direct factor of an open subgroup of $G$.
\end{defn}

Theorem~\ref{intro:wilcomp} shows that just infinite profinite groups have~\propLD.  In fact, we will obtain a characterization of property~\propLD that shows it is closely related to the just infinite property.

The second source of inspiration is the articles \cite{BEW} and \cite{CapDeM} of Barnea--Ershov--Weigel and Caprace--De Medts respectively, which introduce and develop the theory of local isomorphisms of \tdlc groups.  This theory is well-behaved for \tdlc groups $G$ where the \defbold{quasi-centre} $\QZ(G)$, that is, the group of elements with open centralizer in $G$, is trivial.  In that case, there is a group $\ms{L}(G)$, the \defbold{group of germs of $G$}, which contains as an open subgroup every \tdlc group that is locally isomorphic to $G$ and has trivial quasi-centre.  In \cite{BEW}, the normal subgroup structure of $G = \ms{L}(H)$ was considered where $H$ is a \hji profinite group that is not virtually abelian.  Write $\Res(G)$ for the intersection of all open normal subgroups of $G$.

\begin{thm}[{\cite[Proposition~5.1 and Theorem~5.4]{BEW}\footnote{The theorem was originally stated for $G = \ms{L}(H)$, but the argument only uses the more general hypothesis.}}]\label{intro:bew}
Let $H$ be a \hji profinite group that is not virtually abelian and let $G$ be a \tdlc group containing $H$ as an open subgroup.  Then $\QZ(H)=\triv$ and exactly one of the following holds:
\begin{description}
\item[(Residually discrete type)] $\Res(G/\QZ(G)) = \triv$;
\item[(Mysterious type)] $\Res(G/\QZ(G)) \neq \triv$ but $\Res(\Res(G/\QZ(G))) = \triv$;
\item[(Simple type)] $\Res(G/\QZ(G))$ is open and topologically simple.
\end{description}
\end{thm}

In \cite{BEW}, the authors also comment that they are not aware of any groups of germs of \hji profinite groups of mysterious type.  As far as I am aware, it is still unknown whether there are \tdlc groups of mysterious type (not just groups of germs) with trivial quasi-centre and a \hji compact open subgroup. 

In the present article we continue this analysis to other classes of \tdlc groups $G$ that have strong restrictions on their locally normal subgroups, including the case when $G$ is locally isomorphic to a just infinite branch group.

The third source of inspiration is the series of articles \cite{CRW-TitsCore}, \cite{CRW-Part1} and \cite{CRW-Part2} of P.-E. Caprace, G. Willis and the present author, in which further methods were developed for using local properties of \tdlc groups $G$ to obtain restrictions on the global structure (particularly when $G$ is nondiscrete, compactly generated and topologically simple).  Particularly relevant here is the notion of the local decomposition lattice of a \tdlc group $G$, which was also directly inspired by Wilson's approach to just infinite groups.

\begin{defn}
Let $G$ be a \tdlc group.  A subgroup is \defbold{locally normal} if it has open normalizer, and two subgroups $H$ and $K$ are \defbold{locally equivalent} if $H \cap K$ is open in both $H$ and $K$.  The \defbold{structure lattice} $\lnorm(G)$ of $G$ is the bounded lattice formed by the local equivalence classes of locally normal subgroups of $G$, ordered by inclusion of representatives; write $0$ for the least element $[\triv]$ and $\infty$ for the greatest element $[G]$.  When $\QZ(G)$ is discrete, the \defbold{local decomposition lattice} $\ldlat(G)$ consists of those elements $\alpha$ of $\lnorm(G)$ with a complement in $\lnorm(G)$, that is, $\beta:= \alpha^\bot$ such that $\alpha \vee \beta = \infty$ and $\alpha \wedge \beta = 0$.
\end{defn}

Note that for a \tdlc group $G$ such that $\QZ(G)$ is discrete, property~\propLD is exactly the condition that $\lnorm(G) = \ldlat(G)$.

\subsection{Main results}

Our first main result is to characterize property~\propLD, establishing the connection with just infinite profinite groups.

\begin{thm}[See Section~\ref{sec:LD}]\label{thm:propld}Let $G$ be a \tdlc group.  The following are equivalent:
\begin{enumerate}[(i)]
\item $G$ has~\propLD;
\item $G$ locally isomorphic to a profinite group of the form
$$ \prod_{i \in I}L_i,$$
such that finitely many factors $L_i$ (possibly none) are just infinite profinite groups and the remaining factors are finite simple groups.
\end{enumerate}
\end{thm}

We also note that in a group with~\propLD, every closed subnormal subgroup has~\propLD and has an open subgroup that is locally normal: see Lemma~\ref{lem:ld:hered}.

At this point it is useful to introduce a notion of homogeneity of locally normal subgroups, which will capture many examples of groups with~\propLD.

\begin{defn}Let $G_1$ and $G_2$ be \tdlc groups.  Say $G_1$ and $G_2$ are \defbold{commensurable} if there is an isomorphism from a finite index open subgroup of $G_1$ to a finite index open subgroup of $G_2$.  Say $G_1$ and $G_2$ are \defbold{similar} if there are natural numbers $a_1$ and $a_2$ such that $G^{a_1}_1$ is commensurable with $G^{a_2}_2$.  Say $G_1$ and $G_2$ are \defbold{locally similar} if they have compact open subgroups $U_1$ and $U_2$ respectively such that $U_1$ is similar to $U_2$.

A \tdlc group $G$ is \defbold{(locally) monomial} if it is nondiscrete and (locally) similar to every infinite closed locally normal subgroup of $G^a$ for all natural numbers $a$.
\end{defn}

It is easily seen that similarity and local similarity are equivalence relations, stable under passage to an open subgroup.  It is also easy to see that among profinite groups, the monomial property is a similarity invariant.  We can characterize the locally monomial \tdlc groups as follows.  In the next theorem ``LS'' is for ``locally semisimple'', the ``J'' stands for ``just infinite'', and then the letter after indicates the relevant type of just infinite profinite group (either (virtually) \uline{a}belian, (similar to) \uline{h}ereditarily just infinite, or (virtually) \uline{b}ranch).

\begin{thm}[See Section~\ref{sec:similar}]\label{thm:monomial}
Let $G$ be a nondiscrete \tdlc group.  Then the following are equivalent:
\begin{enumerate}[(i)]
\item $G$ is locally monomial and has~\propLD;
\item One of the following holds:
\begin{description}
\item[(LS)] Some open subgroup of $G$ is a direct product of $\aleph_0$ copies of a finite simple group;
\item[(JA)] $G$ has a compact open subgroup isomorphic to $\bZ^d_p$ for some natural number $d$ and prime $p$;
\item[(JH)] $\QZ(G)$ is discrete, $\lnorm(G)$ is finite and an open subgroup of $G$ is a direct product of $n$ copies of a \hji profinite group, where $2^n = |\lnorm(G)|$;
\item[(JB)] $\QZ(G)$ is discrete, $\lnorm(G)$ is infinite, $G/\QZ(G)$ acts faithfully on $\lnorm(G)$ and $G$ is locally isomorphic to a just infinite profinite branch group.
\end{description}
\end{enumerate}
\end{thm}

Given a \tdlc group $G$ with property~\propLD, then there is a compact open subgroup that decomposes into a direct product of monomial factors plus a leftover factor $G_{\infty}$, which is either trivial or a direct product of infinitely many finite simple groups, such that each isomorphism type only appears finitely many times.  Moreover, this factorization is stable under local isomorphisms, and for each similarity type it can be used to canonically associate a locally monomial \tdlc group of that similarity type with $G$.  (See Propositions~\ref{aji:reduced:new} and~\ref{prop:monomial_constituent} for the exact statements.)  For the global structure it is convenient to reduce to the case that $\QZ(G)=\triv$, a reduction that can be accomplished as follows.

\begin{lem}[See Section~\ref{sec:monomial_constituents}]\label{lem:monomial:qz}
Let $G$ be a \tdlc group with property~\propLD.  Then $G$ has a characteristic closed subgroup $Q$ such that $G/Q$ has~\propLD, $\QZ(G/Q) = \triv$, and $\ol{\QZ(Q)}$ is open in $Q$.
\end{lem}

For the rest of this introduction we will assume that $\QZ(G)=\triv$.  In that case we have a locally monomial decomposition as follows.

\begin{thm}[See Section~\ref{sec:monomial_constituents}]\label{thm:monomial:global}
Let $G$ be a \tdlc group with property~\propLD such that $\QZ(G)=\triv$.  Then $G$ is first-countable and has an open subgroup of the form
\[
M = M_1 \times M_2 \times \dots \times M_n,
\]
where for $1 \le i \le n$ the groups $M_i$ have the following properties:
\begin{enumerate}[(i)]
\item $M_i$ is closed and characteristic in every closed locally normal subgroup of $G$ containing $M_i$.
\item $M_i$ is locally isomorphic to a just infinite profinite group $J_i$, so in particular $M_i$ is locally monomial; additionally $M_i$ contains every closed locally normal subgroup of $G$ that is locally similar to $J_i$.
\item In the case that $G$ contains a compact just infinite representative $J_i$ of $[M_i]$, then $J_i \le M_i$ and every nontrivial closed normal subgroup of $M_i$ is open in $M_i$.
\item If $G = \ms{L}(G)$, then $G = M$, $\ms{L}(M_i) = M_i$ for $1 \le i \le n$, and $M_i$ has a just infinite compact open subgroup for all $1 \le i \le n$.
\end{enumerate}
\end{thm}

If $G$ has trivial quasi-centre and is locally isomorphic to a direct product of $n$ \hji profinite groups, then $G$ has an open characteristic subgroup that splits as a direct product of $n$ locally \hji groups, and the factors fall into types as in Theorem~\ref{intro:bew}; see Theorem~\ref{thm:hji:atomic}.  We do not resolve the mysterious type question, but we can put some further restrictions on the mysterious type if it does occur.

\begin{thm}[See Section~\ref{sec:hji}]\label{thm:hji:mysterious}
Let $G$ be a \tdlc group with $\QZ(G)=\triv$, such that $G$ has a \hji compact open subgroup and $G$ is compactly generated.  Suppose that $R = \Res(G)$ is open but $\Res(R)$ is trivial.  Then the following holds.
\begin{enumerate}[(i)]
\item There is a noncompact open normal subgroup $K$ of $R$ such that, given any nontrivial closed normal subgroup $Q$ of $R$, then $Q$ contains a $G$-conjugate of $K$.
\item There is a $G$-conjugate of $K$ that properly contains $K$.
\item Given a compactly generated subgroup $H$ of $G$ such that $R \le H$, then $\Res_R(H)$ is open in $R$.  In particular, $R$ is not compactly generated.
\item Let $L$ be a subgroup of $G$ containing $R$ such that $L/R$ is virtually polycyclic.  Then every open normal subgroup of $R$ contains an open normal subgroup of $L$; in particular, $L$ is residually discrete.
\end{enumerate}
\end{thm}

If $G$ is locally isomorphic to a just infinite profinite branch group, we decompose $\Res(G)$ into noncompact, directly indecomposable parts, each of which is either residually discrete or topologically simple; see Theorem~\ref{thm:locaji:structure}.

We conclude the introduction with two more theorems, this time concerning the relationship between topologically simple locally normal subgroups (which are necessarily nondiscrete and noncompact in the present setting), contraction groups, and the intersection of open subgroups normalized by a given subgroup.

\begin{defn}Let $G$ be a \tdlc group.  The \defbold{contraction group} of $g \in G$ is the group
\[
\con(g) := \{ x \in G \mid g^nxg^{-n} \rightarrow 1 \text{ as } n \rightarrow \infty \}.
\]

Given $H \le G$, we define the \defbold{relative Tits core} $G^\dagger_H$ of $H$ acting on $G$ as
\[
G^\dagger_H = \langle \overline{\con(h)} \mid h \in H \rangle;
\]
we also write $G^\dagger := G^\dagger_G$, and for a single element $g$ we define $G^\dagger_g = G^\dagger_{\grp{g}}$.

Define $\Res_G(H)$ to be the intersection of all open $H$-invariant subgroups of $G$.
\end{defn}

In the following two theorems ``finitely many'' includes the possibility of ``zero'', i.e. the group in question is trivial.  Note that if $H$ has a cocompact polycyclic subgroup, then $\Res_G(H) = \ol{G^\dagger_H}$: see \cite[Corollary~1.12]{ReidFlat}.

\begin{thm}[See Section~\ref{sec:simple}]\label{thm:simpleres}
Let $G$ be a \tdlc group with~\propLD such that $\QZ(G)=\triv$. and let $H$ be a compactly generated subgroup of $G$.  Suppose that at least one of the following holds:
\begin{enumerate}[(i)]
\item $H\Res_G(H)/\Res_G(H)$ is virtually polycyclic;
\item $G$ has no hereditarily just infinite compact locally normal subgroup.
\end{enumerate}
Then $\Res_G(H)$ is a direct product of finitely many topologically simple groups, each of which is locally normal in $G$.
\end{thm}

\begin{thm}[See Section~\ref{sec:simple}]\label{thm:aji:titscore}
Let $G$ be a \tdlc group with~\propLD such that $\QZ(G)=\triv$.  Then $\ol{G^\dagger}$ is a direct factor of $\Res(G)$ and is a direct product of finitely many topologically simple groups.  If $G$ is compactly generated, then $\Res(G)/\ol{G^\dagger} \cong \CC_{\Res(G)}(\ol{G^\dagger})$ is locally isomorphic to a direct product of finitely many \hji profinite groups.
\end{thm}

\paragraph{Acknowledgements}
This article started as an offshoot of the project with Pierre-Emmanuel Caprace and George Willis that led to \cite{CRW-TitsCore}, \cite{CRW-Part1} and \cite{CRW-Part2}, and was also partly developed during research visits to Alejandra Garrido and John Wilson at Oxford, and to Yiftach Barnea at RHUL.  I thank all of them for their hospitality and for very helpful discussions.

\section{Preliminaries}

\subsection{The structure lattice and the decomposition lattice}\label{sec:lnorm}

We begin by setting some terminology.

\begin{defn}
Write $H \le_o G$ to mean ``$H$ is an open subgroup of $G$''.  All groups in this article will be locally compact Hausdorff topological groups.

We distinguish between two groups $H$ and $K$ being \defbold{commensurable}, meaning there is an isomorphism from a finite index open subgroup of $H$ to a finite index open subgroup of $K$, and \defbold{commensurate}, meaning they are both subgroups of some ambient group, such that $H \cap K$ is open and has finite index in both $H$ and $K$.  An automorphism $\alpha$ of some $G \ge H$ \defbold{commensurates} $H$ if $\alpha(H)$ is commensurate with $H$.

Given \tdlc groups $H$ and $K$, a \defbold{local isomorphism} from $H$ to $K$ is a continuous injective open homomorphism $\theta: U \rightarrow K$, where $U \le_o H$.  Two \tdlc groups are \defbold{locally isomorphic} if a local isomorphism between them exists.

Let $G$ be a \tdlc group.  Given two subgroups $H$ and $K$ of $G$, we say $H$ is \defbold{locally equivalent} to $K$, and write $H \sim_o K$, if $H \cap K$ is open in both $H$ and $K$.  Write $[H]$ for the local equivalence class of $H$.  Say $H \le G$ is \defbold{locally normal} if $\N_G(H)$ is open in $G$.  The \defbold{structure lattice} of $G$ is defined to be the set
\[
\lnorm(G) := \{[H] \mid H \le G \text{ closed}, \N_G(H) \le_o G\},
\]
partially ordered by setting $[H] \le [K]$ if $H \le K$.  We write $0$ for the least element $[\triv]$ and $\infty$ for the greatest element $[G]$.

The \defbold{quasi-centre} $\QZ(G)$ of a \tdlc group $G$ is the set of all elements with open centralizer; we say $G$ is \defbold{quasi-discrete} if $\QZ(G)$ is dense in $G$.

More generally, the \defbold{quasi-centralizer} $\QC_G(H)$ of $H \le G$ is the set of all elements of $G$ that centralize an open subgroup of $H$.  Given a local equivalence class $\alpha$ of closed subgroups of $G$, we also define $\QC_G(\alpha)$ to be $\QC_G(H)$ where $H \in \alpha$ (note that the choice of representative is irrelevant).  We also define $\CC^2_G(H) := \CC_G(\CC_G(H))$ and $\QC^2_G(\alpha) := \QC_G(\QC_G(\alpha))$.
\end{defn}

In a profinite group, $\QZ(G)$ is just the union of all finite conjugacy classes of $G$.  Let us note two useful consequences of this fact.

\begin{lem}[{See \cite[Theorems~5.1 and 5.2]{Neumann}}]\label{lem:qz:tf}
Let $G$ be a profinite group.  Then every element of $\QZ(G)$ of finite order is contained in a finite normal subgroup of $G$.  If $\QZ(G)$ is torsion-free, then it is abelian.
\end{lem}

The next lemma is clear from the fact that the conjugation action is continuous; in particular, we see that $\QZ(G)$ contains all discrete locally normal subgroups of a \tdlc group $G$.

\begin{lem}\label{lem:qz:discrete}
Let $G$ be a \tdlc group, let $O \le_o G$ and let $g \in G$.  Suppose that the set $\{ogo\inv \mid o \in O\}$ is discrete.  Then $g \in \QZ(G)$.
\end{lem}

If we rule out quasi-central and abelian locally normal subgroups, we obtain a Boolean algebra canonically associated to the structure lattice.

\begin{defn}
A \tdlc group $G$ is \defbold{[A]-semisimple} if $\QZ(G)=\triv$ and there are no nontrivial abelian locally normal subgroups of $G$.
\end{defn}

\begin{defn}
Given a bounded lattice $\mc{L}$, a \defbold{pseudocomplement} is a map $\perp: \mc{L} \rightarrow \mc{L}$ with the property
\[
\forall \alpha,\beta \in \mc{L}: \alpha \wedge \beta = 0 \Leftrightarrow \beta \le \alpha^\perp.
\]
If it exists, the pseudocomplement is unique, and we call $\mc{L}$ a \defbold{pseudocomplemented lattice}.  A pseudocomplement is a \defbold{complement} if in addition $\alpha \vee \alpha^\perp = \infty$.  In particular, a \defbold{Boolean algebra} is a bounded distributive complemented lattice.
\end{defn}

The set of pseudocomplements of any bounded lattice forms a Boolean algebra: see \cite{Glivenko}.

\begin{thm}[{See \cite[Theorems~3.19 and~5.2]{CRW-Part1}}]\label{thm:cent_pseudocomp}
Let $G$ be an [A]-semisimple \tdlc group.  Then there is a well-defined pseudocomplement $\perp$ on $\lnorm(G)$ given by $[K]^\perp = [\CC_G(K)]$.  Moreover, given closed locally normal subgroups $K$ and $L$ of $G$, then the following are equivalent:
\begin{enumerate}[(i)]
\item $K \cap L$ is discrete, that is, $[K] \wedge [L] = 0$;
\item $K \cap L$ is trivial;
\item $K \le \QC_G(L)$ and $L \le \QC_G(K)$;
\item $[K,L] = \triv$.
\end{enumerate}
\end{thm}

Note that, since condition (i) in Theorem~\ref{thm:cent_pseudocomp} only depends on the local equivalence classes of $K$ and $L$, the same is true for the other three conditions.

The image $\perp\!(\lnorm(G))$ with induced partial order is the \defbold{centralizer lattice} $\lcent(G)$.  The operation $\perp$ is an involution on $\lcent(G)$, that is, $\alpha^{\perp \perp \perp} = \alpha^{\perp}$ for all $\alpha \in \lnorm(G)$, and it serves as the complement for $\lcent(G)$.  The theorem also yields a global version of the centralizer lattice: we have a $G$-equivariant isomorphism of Boolean algebras
\[
\lcent(G) \rightarrow \LCent(G):= \{\CC_G(K) \mid K \text{ is a closed locally normal subgroup of } G\}; \; \alpha \mapsto \QC^2_G(\alpha)),
\]
where $\LCent(G)$ is ordered by inclusion, and the group $\QC^2_G(\alpha) = \CC_G(\QC_G(\alpha))$ in fact contains every closed locally normal subgroup of $G$ that represents $\alpha$.

Within $\lnorm(G)$ we also have the elements $\alpha$ that are complemented in $\lnorm(G)$, that is, such that $\alpha \vee \alpha^{\perp} = \infty$ in $\lnorm(G)$.  These form another Boolean algebra $\ldlat(G)$, the \defbold{local decomposition lattice}, which is a sublattice of $\lnorm(G)$ and a subalgebra of $\lcent(G)$.  One can also define $\LDlat(G)$ as the subalgebra of $\LCent(G)$ corresponding to $\ldlat(G)$; one sees that the elements of $\LDlat(G)$ are all direct factors of open subgroups of $G$, and given $K_1,\dots,K_n \in \LDlat(G)$ such that $K_i \cap K_j = \triv$ for $i \neq j$, then the group generated by the $K_i$ is closed and forms a direct product $K_1 \times \dots \times K_n$.

\subsection{The Mel'nikov subgroup and just infinite groups}

\begin{defn}In this article a direct product always carries the product topology; in particular, the factors are closed subgroups.

A group $H$ is \defbold{finitely decomposable} if $H = H_1 \times \dots \times H_n$ for some $n \ge 1$, such that each $H_i$ cannot be decomposed further as a direct product.  A \tdlc group $G$ is \defbold{locally finitely decomposable} if every compact open subgroup of $G$ is finitely decomposable.

Let $G$ be a profinite group.  The \defbold{Mel'nikov subgroup} $\M(G)$ of $G$ is the intersection of all closed normal subgroups $K$ of $G$ such that $G/K$ is simple.  Say $G$ is \defbold{Mel'nikov-finite} if $G/\M(G)$ is finite.  A \tdlc group $L$ is \defbold{locally Mel'nikov-finite} if every compact open subgroup of $L$ is Mel'nikov-finite.
\end{defn}

We list here some useful facts about the Mel'nikov subgroup of a profinite group.

\begin{prop}\label{melfacts}Let $G$ be a profinite group.
\begin{enumerate}[(i)]
\item We have $\M(G) < G$ unless $G$ is the trivial group.  In particular, if $G$ is topologically characteristically simple, then $\M(G)=\triv$.
\item Given a closed normal subgroup $H$ of $G$, we have $\M(G/H)=\M(G)H/H$.  In particular, $\M(G/\M(G))=\triv$.
\item The quotient $G/\M(G)$ is a direct product of finite simple groups.  In particular, the torsion elements of $\QZ(G/\M(G))$ form a dense subgroup of $G/\M(G)$.
\item If $G/\M(G)$ is finite, then $G$ is finitely decomposable.
\item Let $H$ be a closed normal subgroup of $G$.  Then $G = \M(G)H$ if and only if $G = H$.
\item Let $H$ be a closed normal subgroup of $G$.  Then $\M(H) \le \M(G)$.  In particular, it follows that $G$ is locally Mel'nikov-finite if and only if there is a set of open normal Mel'nikov-finite subgroups of $G$ forming a basis of identity neighbourhoods.
\item Define $M_0 = G$ and thereafter $M_{i+1} = \M(M_i)$.  Then $\bigcap_{i \in \bN}M_i = \triv$.  Thus if $G$ is locally Mel'nikov-finite, then it is first-countable, and indeed has a countable descending chain of open characteristic subgroups with trivial intersection.
\item Suppose that $G$ is locally Mel'nikov-finite.  Then for each natural number $n$, there are only finitely many open subgroups of $G$ of index $n$.
\end{enumerate}
\end{prop}

\begin{proof}For (iii), see for instance \cite[Lemma~8.2.2]{RZ}.  The other statements are straightforward exercises given the well-known properties of profinite groups.\end{proof}

\begin{defn}
A profinite group $G$ is \defbold{just infinite} if it is infinite, but every nontrivial closed normal subgroup of $G$ is open (in other words, has finite index).  It is \defbold{hereditarily just infinite} (\hji) if every open subgroup of $G$ is just infinite.

A profinite group $G$ is \defbold{virtually just infinite} if $G$ has an open subgroup in common with a just infinite profinite group.  (Virtually) hereditarily just infinite groups and (virtually) just infinite branch groups are defined similarly.
\end{defn}

There are very few possibilities for the quasi-centre of a just infinite profinite group.

\begin{lem}\label{ji:qz}Let $G$ be a just infinite profinite group.  The quasi-centre $\QZ(G)$ is the unique largest abelian closed locally normal subgroup of $G$; in particular, $\QZ(G)$ is closed.  If $\QZ(G)$ is nontrivial then it is a direct product of finitely many copies of $\bZ_p$ for some fixed prime $p$.\end{lem}

\begin{proof}Since $G$ has no finite normal subgroups, Lemma~\ref{lem:qz:tf} implies that $\QZ(G)$ is torsion-free abelian.  The closure $Q$ of $\QZ(G)$ is an abelian closed normal subgroup, which is then either trivial or of finite index; in either case, $Q \le \QZ(G)$, so $\QZ(G)$ is closed.  From now on we may assume $Q > \triv$.  Using the just infinite property, it is now straightforward to see that $Q \cong (\bZ_p)^n$ for some prime $p$ and natural number $n$.

Given an abelian locally normal subgroup $B$ of $G$, we claim that $B$ commutes with $R = \N_Q(B)$, so that $B \le Q$.  Letting $x \in B$ and $k \in \bN$, we see that $[x^k,r] = [x,r]^k$ using the fact that $[B,R]$ is central in $B$.  For some $k > 0$ we have $x^k \in R$, so $[x^k,r]=1$ and hence $[x,r]^k = 1$.  Since $R$ is torsion-free it follows that $[x,r]=1$ as required.
\end{proof}

\begin{cor}\label{cor:ji:Cstable}Let $G$ be a just infinite profinite group.  Then either $G$ is virtually abelian or $G$ is [A]-semisimple.\end{cor}

\subsection{Contraction in \tdlc groups}

We recall some definitions relating to contraction groups in \tdlc groups.

\begin{defn}The \defbold{contraction group} of an automorphism $\phi$ of a group $G$ is the group
\[
\con(\phi) := \{ x \in G \mid \phi^n(x) \rightarrow 1 \text{ as } n \rightarrow \infty \}.
\]

Given $g \in G$, we let $g$ get on $G$ by left conjugation.  
Let $G$ be a \tdlc group and let $H \le \Aut(G)$.  We define the \defbold{relative Tits core} $G^\dagger_H$ of $H$ acting on $G$ as
\[
G^\dagger_H = \langle \overline{\con(h)} \mid h \in H \rangle.
\]

Given a closed subgroup $K$ of $G$ and $H \le \N_G(K)$, we define $K^\dagger_H$ to be $K^\dagger_{\phi(H)}$, where $\phi$ is the left conjugation action of $H$ on $K$.  In particular, we define the \defbold{Tits core} $G^\dagger := G^\dagger_G$.  For a single automorphism $\alpha$, we write $G^\dagger_{\alpha} := G^\dagger_{\grp{\alpha}}$.
\end{defn}

The Tits core is not sensitive to passing to an open subgroup of finite index, as follows from the following basic observation.

\begin{lem}\label{lem:contraction:powers}Let $G$ be a \tdlc group, let $\alpha \in \Aut(G)$ and let $n$ be a positive integer.  Then
\[
\con(\alpha) = \con(\alpha^n).
\]
\end{lem}

\begin{proof}We have $\con(\alpha^n) \ge \con(\alpha)$, since $(\alpha^{ni})_{i \in \bN}$ is a subsequence of $(\alpha^i)_{i \in \bN}$.  Let $H = \overline{\con(\alpha)}$, let $x \in \con(\alpha^n)$ and set $x_i = \alpha^i(x)$.  Then the sequence $(x_{ni})_{i \in \bN}$ converges to the identity.  Since $\alpha$ is a continuous automorphism and $\alpha^j(x_k) = x_{j+k}$ for all $j,k \in \bZ$, it follows that the sequence $(x_{j+ni})_{i \in \bN}$ converges to $\alpha^j(1) = 1$.  Hence $(x_i)_{i \in \bN}$ converges to the identity, since it can be partitioned into finitely many subsequences $(x_{j+ni})_{i \in \bN}$ for $0 \le j < n$, each of which converges to the identity.  In other words, $x \in \con(\alpha)$.
\end{proof}

For individual automorphisms, having trivial relative Tits core is equivalent to having small invariant neighbourhoods, as follows.

\begin{lem}\label{lem:anisotropic}
Let $G$ be a \tdlc group and let $\alpha$ be an automorphism of $G$.  Then $G^\dagger_{\alpha} = \triv$ if and only if the compact open subgroups $U$ of $G$ such that $\alpha(U) = U$ form a base of neighbourhoods of the identity.
\end{lem}

\begin{proof}
If $G^\dagger_{\alpha} = \triv$, then by \cite[Proposition~3.24]{BW}, $\alpha$ normalizes a compact open subgroup $U$, and by \cite[Theorem~3.32]{BW}, $U$ can be made arbitrarily small.  Conversely, if the compact open subgroups $U$ of $G$ such that $\alpha(U) = U$ form a base of neighbourhoods of the identity, then clearly $\con(\alpha^n) = \triv$ for all $n \in \bZ$, so $G^\dagger_{\alpha} = \triv$.
\end{proof}

Given a group $G$ and a subgroup $H$, write $\Res_G(H)$ for the intersection of all open $H$-invariant subgroups of $H$.

\begin{lem}[{See \cite[Theorem~B]{ReidDistal}}]\label{lem:reduced_envelope}
Let $G$ be a \tdlc group and let $H$ be a compactly generated subgroup of $G$.  Then there is a compactly generated open subgroup $E$ of $G$ such that $H \le E$ and $\Res_G(H) = \Res(E)$.
\end{lem}

\section{Groups in which every locally normal subgroup is locally a direct factor}\label{sec:LD}

\subsection{First observations}

The goal of this section is to prove Theorem~\ref{thm:propld}.  We start by establishing some general features of the property~\propLD.

\begin{lem}\label{lem:ld:hered}
Let $G$ be a \tdlc group with~\propLD and let $K$ be a closed subnormal subgroup of an open subgroup of $G$.  Then some open subgroup of $K$ is a local direct factor of $G$; moreover, $K$ has~\propLD.
\end{lem}

\begin{proof}
Let $U \le_o G$ such that $K$ is subnormal in $U$ and let
$$ K = U_0 \unlhd U_1 \unlhd \dots \unlhd U_n = U$$
be a subnormal series from $K$ to $U$ of shortest possible length.  By induction on $n$, there is an open subgroup $V$ of $U_1$ that is a local direct factor of $U$, so there is $M \le U$ such that $\langle V,M \rangle \cong V \times M \le_o U$.  Now $K \cap V$ is normal in $\langle V,M \rangle$, hence locally normal in $G$, so $G$ has an open subgroup $K_0 \times N$ where $K_0 \le_o K \cap V$.  In particular, $K_0$ is an open subgroup of $K$ that is a local direct factor of $G$.

Every locally normal subgroup $R$ of $K_0$ is also locally normal in $G$, and so there is $R_2 \le_o R$ and $L_2 \le G$ such that $R_2 \times L_2 \le_o G$.  We now see that $R_2 \times (L_2 \cap K_0) \le_o L$, so $R_2$ is a local direct factor of $K_0$.  Thus $K_0$ has \propLD; since $K_0$ is open in $K$, it follows that $K$ has \propLD.
\end{proof}

When the quasi-centre is trivial, property~\propLD can be characterized in terms of quasi-centralizers.

\begin{lem}\label{lem:ld:qz}
Let $G$ be a \tdlc group with~\propLD and let $K$ be a closed locally normal subgroup of $G$ such that $\QZ(G) \cap K = \triv$.  Then $\QZ(K) = \triv$.
\end{lem}

\begin{proof}
Without loss of generality, $K$ is normal in $G$.  Let $K_0$ be an open subgroup of $K$ that is a direct factor of an open subgroup of $G$.  Then every element of $\QZ(K_0)$ has open centralizer in $G$, so we have
\[
\QZ(K) \cap K_0 = \QZ(K_0) \le \QZ(G) \cap K = \triv.
\]
Thus $\QZ(K)$ is a discrete normal subgroup of $G$; hence $\QZ(K) \le \QZ(G)$ by Lemma~\ref{lem:qz:discrete}, so in fact $\QZ(K) = \triv$.
\end{proof}

\begin{lem}\label{lcent:ldlat}Let $G$ be a \tdlc group with $\QZ(G)=\triv$.  The following are equivalent:
\begin{enumerate}[(i)]
\item $G$ has~\propLD;
\item $G$ is [A]-semisimple and for every compact locally normal subgroup $K$ of $G$ that is not open, then $\QC_G(K)$ is nondiscrete.
\end{enumerate}
\end{lem}

\begin{proof}
Let $K$ be a compact locally normal subgroup of $G$ that is not open.

Suppose $G$ has~\propLD.  By Lemma~\ref{lem:ld:qz}, every closed locally normal subgroup of $G$ has trivial quasi-centre, so $G$ is [A]-semisimple.  After replacing $K$ with an open subgroup we find that $K \times \QC_G(K)$ is open in $G$, so certainly $\QC_G(K)$ is not discrete.  Thus (ii) holds.

Conversely, suppose (ii) holds.  Given a compact open subgroup $U$ containing $K$, we have a compact locally normal subgroup of the form $L = K \times \QC_U(K)$.  Because $[\QC_U(K)]$ is the pseudocomplement of $[K]$ in $\lnorm(G)$ (see Section~\ref{sec:lnorm}), we see that $[\QC_U(L)] = 0$, so by (ii), $L$ must be open.  Thus $K$ is a direct factor of an open subgroup of $G$, showing that $G$ has~\propLD.
\end{proof}

We have the following restriction on local direct factors.

\begin{lem}\label{lem:melfin}Let $G$ be a \tdlc group with~\propLD such that $\QZ(G)=\triv$.  Then $K/\overline{[K,K]}$ and $K/\M(K)$ are both finite for every compact locally normal subgroup $K$ of $G$.  In particular, $G$ is locally finitely decomposable and first-countable.\end{lem}

\begin{proof}
Let $K$ be a compact locally normal subgroup of $G$.   Then $K$ has \propLD by Lemma~\ref{lem:ld:hered} and trivial quasi-centre by Lemma~\ref{lem:ld:qz}.

By property~\propLD, $K$ has an open subgroup of the form $L_1 \times L_2$ where $L_1 \le _o \overline{[K,K]}$.  Since $L_1$ has finite index in $\overline{[K,K]}$, $L_2$ is virtually abelian, and indeed $L_2$ can be chosen to be abelian.  But then $L_2 \le \QZ(K) = \triv$, so in fact $L_1$ is open in $K$, and hence $\overline{[K,K]}$ is open in $K$.  Similarly, $K$ has an open subgroup of the form $M_1 \times M_2$ where $M_1 \le_o \M(K)$.  We see that $M_2$ is commensurable with a profinite group $R$ such that $\M(R)=\triv$.  Then $\QZ(R)$ is dense in $R$ by Proposition~\ref{melfacts}(iii).  But then $\overline{\QZ(M_2)}$ must be open in $M_2$; since $\QZ(M_2) = \triv$, we conclude that $M_2$ is finite and $M_1$ is open in $K$, so $\M(K)$ is open in $K$.

Applying parts (iv) and (vii) of Proposition~\ref{melfacts} to a compact open subgroup, we see that $G$ is locally finitely decomposable and first-countable.
\end{proof}

The next lemma will allow us to split the proof of Theorem~\ref{thm:propld} into the quasi-discrete case and the trivial quasi-centre case.

\begin{lem}\label{lem:propd:decomp}Let $G$ be a \tdlc group.
\begin{enumerate}[(i)]
\item If $G$ has~\propLD, then $G$ has a compact open subgroup of the form $Q \times R$ where $Q = \ol{\QZ(Q)}$ and $\QZ(R)=\triv$, such that $Q$ and $R$ both have~\propLD.
\item If $G$ has a compact open subgroup $Q \times R$ where $Q$ and $R$ have~\propLD and $\QZ(R)=\triv$, then $G$ has~\propLD.
\end{enumerate}
\end{lem}

\begin{proof}
Suppose $G$ has~\propLD and let $K = \ol{\QZ(G)}$.  Then $G$ has a compact open subgroup of the form $Q \times R$ where $Q$ is open in $K$.  We immediately see that $Q = \ol{\QZ(Q)}$.  Since $Q$ already accounts for an open subgroup of $K$, the intersection $K \cap R$ is discrete, so $\QZ(G) \cap R$ is discrete, and by passing to an open subgroup of $R$ we may ensure that $\QZ(G) \cap R = \triv$, so $\QZ(R)=\triv$ by Lemma~\ref{lem:ld:qz}.  By Lemma~\ref{lem:ld:hered}, $Q$ and $R$ have~\propLD.

Now suppose the hypothesis of (ii) holds.  Since property~\propLD is a local property we may assume $G = Q \times R$.  Let $K$ be a closed locally normal subgroup of $G$; we must show that $K$ is locally equivalent to a local direct factor of $G$.  Let $K_Q$ and $K_R$ be the projections of $K$ onto $Q$ and $R$ respectively.  By passing to an open subgroup, we may assume $K_Q$ is a direct factor of $Q$ and $K_R$ is a direct factor of $R$.  Note also that $K_Q$ and $K_R$ have \propLD by Lemma~\ref{lem:ld:hered}, and $\QZ(K_R)=\triv$ by Lemma~\ref{lem:ld:qz}.  So it suffices to consider the case when $K_R = R$.

By Lemma~\ref{lem:melfin}, $\ol{[R,R]}$ is open in $R$.  Moreover, since $R \le KQ$ and $K$ is normal in $G$, we have $\ol{[R,R]} = \ol{[K,R]} \le K$, so indeed $K \cap R$ is open in $R$.  It follows that the subgroup $(K \cap Q) \times (K \cap R)$ is open in $K$.  Using property~\propLD in $Q$, we obtain an open subgroup of $G$ of the form $Q_2 \times K_2 \times (K \cap R)$ where $Q_2 \le Q$ and $K_2 \le_o K \cap Q$.  Thus $G$ has \propLD as required.
\end{proof}

\subsection{[A]-semisimple groups with property~\propLD}

Lemma~\ref{lem:melfin} is a strong restriction on locally normal subgroups that will lead to a characterization of property~\propLD in the [A]-semisimple case.

\begin{lem}\label{lem:quasitrans}
Let $G$ be an [A]-semisimple profinite group such that every closed normal subgroup is finitely decomposable.  Then $\ldlat(G)^G$ is finite.
\end{lem}

\begin{proof}
Suppose $\ldlat(G)^G$ is infinite and let $U$ be an open normal subgroup of $G$ such that $\QZ(U)=\triv$.  Then we can find an infinite set $\mc{L}$ of pairwise disjoint elements of $\ldlat(G)^G$.  The corresponding direct product
\[
\prod_{\alpha \in \mc{L}}\QC^2_U(\alpha)
\]
is then a closed normal subgroup of $G$ that is not finitely decomposable.
\end{proof}

\begin{lem}\label{lem:quasimin}
Let $G$ be an [A]-semisimple profinite group.  Suppose that $N$ is a closed normal subgroup of $G$, such that $[N]$ is a minimal nonzero element of $\lnorm(G)^G$ and such that $N\CC_G(N)$ is open in $G$.  Then $N$ has~\propLD and is commensurable with a just infinite profinite group.
\end{lem}

\begin{proof}
To show $N$ has \propLD, by Lemma~\ref{lcent:ldlat} it suffices to show that every closed locally normal subgroup of $N$ of infinite index has nondiscrete (quasi-)centralizer in $N$.  Let $K$ be a closed locally normal subgroup of $N$ of infinite index.  By replacing $N$ and $K$ with $N \cap U$ and $K \cap U$ respectively for a sufficiently small open normal subgroup $U$ of $G$, we may assume $K$ is normal in $U$ and $N \le U$.  Let $K_1,\dots,K_n$ be the conjugates of $K$ in $G$, with $K_1 = K$.  By the minimality of $[N]$, we have $\bigcap^n_{i=1}K_i = \triv$.  Let $I$ be a subset of $\{1,\dots,n\}$ of largest possible size such that $L := \bigcap_{i \in I}K_i > \triv$; since $I$ is a proper subset of $\{1,\dots,n\}$, we may assume $1 \not\in I$.  Then $L$ is a normal subgroup of $N$ and $K \cap L = \triv$, so $\QC_N(K) \ge L > \triv$.  Indeed $\QC_N(K)$ must be nondiscrete, since $\overline{\QC_N(K)}$ is a nontrivial locally normal subgroup of $G$.  Thus $N$ has \propLD as required.

Since $N\CC_G(N)$ is open in $G$, we can identify $\lnorm(N)$ with the ideal $\mc{I}$ generated by $[N]$ in $\lnorm(G)$; we also note that $N \cap \CC_G(N) \le \QZ(G) = \triv$.  Let $H = G/\CC_G(N)$; we see that there is a $G$-equivariant isomorphism between $\lnorm(H)$ and $\mc{I}$, and $H$ is commensurable with $N$.  If $x \in G$ is such that $x\CC_G(N)$ centralizes an open subgroup of $G/\CC_G(N)$, then $x \in \QC_G(N) = \CC_G(N)$, so $\QZ(H) =\triv$; thus $H$ has no nontrivial finite normal subgroups.  Moreover, given the minimality of $[N]$ and the way $H$ arises as a quotient of $G$, the action of $H$ on $\lnorm(H)$ has no nontrivial fixed points.  Thus $H$ is just infinite.
\end{proof}

\begin{cor}\label{cor:quasimin}
Let $G$ be an [A]-semisimple profinite group with~\propLD.  Then $\lnorm(G)^G$ is finite.  Moreover, there is an open characteristic subgroup $U$ of $G$ that decomposes as a finite direct product, where each factor is commensurable with a just infinite group.
\end{cor}

\begin{proof}
After passing to an open subgroup, we may assume $\QZ(G)=\triv$.  Given a closed locally normal subgroup $H$ of $G$, then by Lemmas~\ref{lem:ld:hered} and~\ref{lem:ld:qz}, $H$ has~\propLD with $\QZ(H) =\triv$, so $\lnorm(H) = \ldlat(H)$ and $H$ is finitely decomposable by Lemma~\ref{lem:melfin}.  It follows by Lemma~\ref{lem:quasitrans} and property~\propLD that $\lnorm(G)^G$ is finite.

Let $\alpha_1,\dots,\alpha_n$ be the atoms of $\lnorm(G)^G$.  By property~\propLD, $\lnorm(G) = \ldlat(G)$, and hence $G$ has an open characteristic subgroup $U = \prod^n_{i=1}H_i$ where $H_i = \QC^2_G(\alpha_i)$.  For $1 \le i \le n$, we see that $H_i$ is commensurable with a just infinite group by Lemma~\ref{lem:quasimin}.
\end{proof}

\subsection{The quasi-discrete case}

Let us now focus on the case of a profinite group $G$ such that $\QZ(G)$ is dense in $G$.  Write $\bP$ for the set of all prime numbers.

\begin{lem}\label{lem:elab:complement}Let $G$ be a profinite group that is elementary abelian, that is, $G$ is abelian and there is some $p \in \bP$ such that $x^p=1$ for all $x \in G$.  Then every closed subgroup of $G$ is a direct factor, with a closed direct complement.  In particular, $G$ has~\propLD.\end{lem}

\begin{proof}The Pontryagin dual $\widehat{G} := \Hom(G,\bT)$ is an elementary abelian discrete group, which we may regard as a vector space over the field of $p$ elements.  Let $K$ be a closed subgroup of $G$.  Then the annihilator $M$ of $K$ in $\widehat{G}$ is a subspace of $\widehat{G}$, so we have $\widehat{G} = M \oplus N$ for some other subspace $N$ of $\widehat{G}$ (for instance, we could take a basis $A$ of $M$, extend to a basis $B$ of $\widehat{G}$, and set $N$ to be the span of $B \setminus A$).  Now the annihilator $R$ of $N$ in $G$, that is, the set of elements $g \in G$ such that $g \in \ker\phi$ for all $\phi \in \widehat{G}$, is a closed subgroup of $G$, and we have $G = K \times R$.\end{proof}

We can now prove Theorem~\ref{thm:propld} for profinite groups with dense quasi-centre.

\begin{prop}\label{prop:ld:qzdense}Let $G$ be a quasi-discrete profinite group.  The following are equivalent:
\begin{enumerate}[(i)]
\item $G$ has~\propLD.
\item $G$ has an open subgroup of the form 
$$ \prod^k_{i=1} (\bZ_{p_i})^{n_i} \times E,$$
where $k \ge 0$, $\{p_1,\dots,p_k\}$ is a set of primes, with associated positive integers $(n_1,\dots,n_k)$ and $E$ is a direct product of finite simple groups.
\end{enumerate}
\end{prop}

\begin{proof}
Suppose that $G$ has~\propLD.  Let $K$ be the closed subgroup of $G$ topologically generated by all minimal finite normal subgroups of $G$.  Every minimal finite normal subgroup of $G$ is a product of copies of a finite simple group; we see from this that $K/K_0$ is a direct product of simple groups for every open normal subgroup $K_0$ of $K$, and hence $\M(K) = \triv$.  Thus by Proposition~\ref{melfacts}(iii), $K = \prod_{i \in I}S_i$ for some finite simple groups $S_i$.  By property~\propLD, there is an open subgroup of $G$ of the form $G_2 = K_2 \times L$ where $K_2 \le_o K$ and $L$ intersects $K$ trivially; we see that $\QZ(G_2/K_2)$ is dense in $G_2/K_2$, so $\QZ(L)$ is dense in $L$.  We can also take $K_2 = \prod_{i \in I'}S_i$, where $I'$ is a cofinite subset of $I$.

Let $x \in \QZ(L)$ and suppose $x$ has finite order.  Then $x \in \QZ(G)$, so by Lemma~\ref{lem:qz:tf}, $x$ is contained in a finite $G$-invariant subgroup $F$ of $L$.  Now $F$ is a finite normal subgroup of $N$ that intersects trivially with $K$.  From the definition of $K$, this forces $F = \triv$, so $x=1$.  We have thus shown that $\QZ(L)$ is torsion-free.  By Lemma~\ref{lem:qz:tf}, $\QZ(L)$ is abelian, hence $L = \overline{\QZ(L)}$ is abelian; in particular, $L = \QZ(L)$.  Now $L$ has \propLD by Lemma~\ref{lem:ld:hered}.  There is therefore an open subgroup of $L$ of the form $L_1 \times L_2$, where $L_1 \le_o \M(L)$.  We see that $L/\M(L)$ is a product of elementary abelian groups, so by passing to an open subgroup, we can make $L_2$ a product of elementary abelian groups.  But $L$ is torsion-free, so then $L_2 = \triv$.  We conclude that $|L:\M(L)| < \infty$; since $L$ is torsion-free abelian, it follows that
\[
L =  \prod^k_{i=1}(\bZ_{p_i})^{n_i}.
\]
for some distinct primes $p_1,\dots,p_k$ ($k \ge 0$) and positive integers $n_1,\dots,n_k$.  Thus (ii) holds.

Conversely, suppose that (ii) holds.  Without loss of generality, let us suppose that
$$ G = T \times A \times E,$$
where $T = \prod^k_{i=1} (\bZ_{p_i})^{n_i}$, $A = \prod_{p \in \bP} A_p$ with $A_p$ elementary abelian and $E$ is a direct product of nonabelian finite simple groups.  Let $K$ be a closed locally normal subgroup of $G$.  On passing to an open subgroup of $G$, we may in fact assume that $K$ is normal in $G$.  In particular, $K_E = K \cap E$ is normal in $E$, from which we may conclude that $\Z(K_E)=\triv$ and also $K = \CC_K(K_E)K_E$ (since $G$ only induces inner automorphisms on the normal subgroups of $E$ by conjugation), so in fact $K = \CC_K(K_E)K_E$.  So to show that $G$ has~\propLD, it suffices to show that $T \times A$ has~\propLD, and thus we may assume $E=\triv$.  In this case, $K$ is abelian, so admits a canonical direct decomposition $K := \prod_{p \in \bP}K_p$ by Sylow's theorem, where $K_p$ is the $p$-Sylow subgroup of $K$.  Note that $T$ is a pro-$\pi$ group for a finite set of primes $\pi$.  We can construct a subgroup $L = \prod_{p \in \bP}L_p$ that is almost a direct complement of $K$ in $G$ (in the sense that an open subgroup of $G$ is of the form $K_2 \times L$, with $K_2 \le_o K$) as follows:

For all primes $p \not\in \pi$, we have $K_p \le A_p$, and we may choose $L_p$ to be a direct complement of $K_p$ in $A_p$ using Lemma~\ref{lem:elab:complement}.

Given $p \in \pi$, choose a direct complement $R_p$ to $K_p \cap A = K_p \cap A_p$ in $A_p$, noting that $K_p \cap A$ is the torsion subgroup of $K_p$ (since $A$ contains the torsion subgroup of $G$).  Let $P$ be the $p$-Sylow subgroup of $G$.  Then $P/A_p \cong (\bZ_p)^n$ for some integer $n$.  The group $(\bZ_p)^n$ has \propLD for the following reason: given a subgroup $P$ of $\bZ^n_p$, we embed $\bZ^n_p$ in the vector space $\bQ^n_p$, obtain a $\bQ_p$-linear complement $Q$ to the span of $P$, and observe that $P \times (Q \cap \bZ^n_p)$ spans $\bQ^n_p$ and is therefore an open subgroup of $\bZ^n_p$.  Thus $P/A_p$ has \propLD.  Note that there is a continuous injective homomorphism from $P/A$ to $P \cap T$ given by $\phi: x \mapsto x^p$.  Let $Q/A_p$ be a closed subgroup such that $Q/A_p \times K_pA_p/A_p$ is open in $P/A_p$, and let $S_p = \phi(Q)$.  Finally, set $L_p = R_p \times S_p$.  We see that $P$ has an open subgroup of the form
$$S_p \times \phi(K_p) \times R_p \times (K_p \cap A_p),$$
and $\phi(K_p) \times (K_p \cap A_p)$ is an open subgroup of $K_p$.

Finally, set $K_2 = \prod_{p \in \pi} (\phi(K_p) \times (K_p \cap A_p)) \times \prod_{p \not\in \pi} K_p$ and $L = \prod_{p \in \bP}L_p$.  Then $K_2 \times L$ is an open subgroup of $G$ with $K_2$ open in $K$, as required.
\end{proof}

\subsection{Characterization of property~\propLD}

We now complete the characterization of property~\propLD.

\begin{proof}[Proof of Theorem~\ref{thm:propld}]
Suppose $G$ has~\propLD.  Then by Lemma~\ref{lem:propd:decomp}, $G$ has a compact open subgroup of the form $Q \times R$ where $Q$ is quasi-discrete and $\QZ(R)=\triv$.  By Proposition~\ref{prop:ld:qzdense}, $Q$ has an open subgroup that is a direct product of finite simple groups (possibly including cyclic groups of prime order) together with finitely many groups $\bZ_{p_i}$, the latter being just infinite profinite groups.  By Corollary~\ref{cor:quasimin}, $R$ has an open subgroup that is commensurable with a finite direct product of just infinite groups.  Thus $G$ has an open subgroup in common with a direct product 
$$ \prod_{i \in I}L_i,$$
such that finitely many factors $L_i$ (possibly none) are just infinite profinite groups and the remaining factors are finite simple groups, as required.

Conversely, suppose $G$ has an open subgroup in common with a direct product 
$$ \prod_{i \in I}L_i,$$
such that finitely many factors $L_i$ (possibly none) are just infinite profinite groups and the remaining factors are finite simple groups.  Without loss of generality we may assume actually $G =  \prod_{i \in I}L_i$.  Recalling from Lemma~\ref{ji:qz} the structure of just infinite groups with nontrivial quasi-centre, we can write an open subgroup of $G$ as a direct product
\[
\prod^k_{i=1}(\bZ_{p_i})^{n_i} \times E \times J_1 \times \dots \times J_n,
\]
where $E$ is a direct product of finite simple groups and each of the groups $J_i$ is just infinite with trivial quasi-centre.  Appealing to Proposition~\ref{prop:ld:qzdense}, the group $Q = \prod^k_{i=1}(\bZ_{p_i})^{n_i} \times E$ has \propLD, whilst $J_i$ has \propLD for all $i$ by Lemma~\ref{lem:quasimin}.  Using Lemma~\ref{lem:propd:decomp}, we conclude that $G$ has \propLD, completing the proof of the theorem.
\end{proof}

\section{Global structure}

We now turn to the study of \tdlc groups with~\propLD, using what we know about the compact open subgroups to deduce global properties.

\subsection{Similarity classes}\label{sec:similar}

We now prove Theorem~\ref{thm:monomial}, starting with two lemmas.

\begin{lem}\label{lem:monomial:ji}
Every virtually just infinite profinite group is monomial.
\end{lem}

\begin{proof}
Since the monomial property is a commensurability invariant, we only need to consider just infinite profinite groups.  Let $J$ be a just infinite profinite group and let $H$ be an infinite closed locally normal subgroup of $J^a$ for some $a \ge 1$; we must show that $H$ is similar to $J$.  If $\QZ(J)>\triv$, then by Lemma~\ref{ji:qz}, $J$ is commensurable with $\bZ^m_p$ for some prime $p$ and natural number $m$, so and then $H$ is commensurable with $\bZ^n_p$ for some $1 \le n \le am$; hence $H$ is similar to $J$.  Thus we may assume $\QZ(J)=\triv$.  In this case, $J^a$ is a finite index subgroup of the just infinite group $J' = J \wr \Sym(a)$, which contains $H$ as a locally normal subgroup. Let $U$ be an open normal subgroup of $J'$ that normalizes $H$; without loss of generality we may replace $H$ with $U \cap H$ and assume that $H$ is normal in $U$.  Then $\lnorm(U)^U$ is finite by Corollary~\ref{cor:quasimin}, so it is generated by its atoms.  We can then take a closed normal subgroup $K$ of $U$ representing an atom of $\lnorm(U)^U$ such that the distinct $J'$-conjugates $K_1,\dots,K_n$ of $K$ represent distinct atoms of $\lnorm(U)^U$.  Since $J'$ is just infinite, $\lnorm(U)^{J'} = \{0,\infty\}$, and so $U$ has an open subgroup of the form $K_1 \times \dots \times K_n$; we then see that every closed normal subgroup of $U$ is commensurate with the product of some subset of $K_1,\dots,K_n$ and hence is similar to $K$.  In particular, $H$ is similar to $K$; clearly also $K$ is similar to $J'$ and $J'$ is similar to $J$.
\end{proof}

\begin{lem}\label{lem:monomial:ji:bis}
Let $G = J_1 \times \dots \times J_n$, where $n \ge 1$ and $J_1,\dots,J_n$ are virtually just infinite profinite groups.  Then the following are equivalent:
\begin{enumerate}[(i)]
\item $G$ is monomial;
\item $J_1,\dots,J_n$ are similar;
\item $G$ is virtually just infinite.
\end{enumerate}
\end{lem}

\begin{proof}
Given Lemma~\ref{lem:monomial:ji}, the only nontrivial implication is that (ii) implies (iii), so let us suppose (ii) holds.  If $\QZ(J_i) > \triv$ for some $i$, then $\QZ(J_j) > \triv$ for $1 \le j \le n$, and hence by Lemma~\ref{ji:qz} and similarity, $G$ is commensurable with $\bZ^m_p$ (for some prime $p$ and natural number $m$), which in turn is commensurable with a just infinite virtually abelian group of the form $\bZ^m_p \rtimes F$ for some finite subgroup $F$ of $\mathrm{GL}_m(\bZ_p)$ acting irreducibly over $\bQ_p$, implying (iii).  Thus we may assume $\QZ(J_i) = \triv$ for $1 \le i \le n$; hence also $\QZ(G) = \triv$.

There is a profinite group $H$ and natural numbers $k_1,\dots,k_n$ such that $J^{k_i}_i$ is commensurable with $H$ for $1 \le i \le n$.  We can clearly take $H$ to be just infinite, for instance $H = J_1 \wr \Sym(k_1)$.  Taking $k = \prod^n_{i=1}k_i$, we see that $G^k$ is commensurable with a power $H^l$ of $H$.  In turn, $H^l$ is commensurable with the just infinite profinite group $L = H \wr \Sym(l)$.

We now argue that $G$ itself is commensurable with a just infinite profinite group.  We have an isomorphism $\theta$ from $V^k$, where $V$ is open in $G$, to an open subgroup $U$ of $L$.  In particular, $\theta(V)$ represents an element of $\lnorm(U)^U$.  Similar to the proof of Lemma~\ref{lem:monomial:ji}, we find that $V$ is commensurable with a direct power $K^m$ of $K$, where $K$ represents an atom of $\lnorm(U)^U$.  It follows by Lemma~\ref{lem:quasimin} that $K$ is commensurable with a just infinite profinite group $K_2$, and thus $G$ is commensurable with the just infinite profinite group $K_2 \wr \Sym(m)$.  Thus (ii) implies (iii), completing the cycle of implications.
\end{proof}

\begin{proof}[Proof of Theorem~\ref{thm:monomial}]
Let us first consider the locally normal subgroups of an infinite direct product $E = \prod_{i \in I} S_i$, where $S_i$ is a finite simple group.  We see that every closed locally normal subgroup of $E$ is commensurate with the subgroup $\prod_{i \in I'}S_i$ for some subset $I'$ of $I$, where $I'$ is determined up to adding or removing finitely many elements.  From this description, we see that if $E$ is monomial, then all but finitely many $S_i$ belong to a single isomorphism class; moreover, the fact that $E$ is similar to $\prod_{i \in I'}S_i$ where $|I'| = \aleph_0$ ensures that $|I| = \aleph_0$.  Conversely, if all but finitely many $S_i$ belong to a single isomorphism class and $I$ is countably infinite, it is clear that $E$ is monomial. 

Suppose $G$ is locally monomial and has~\propLD.  Then given Theorem~\ref{thm:propld}, $G$ is locally isomorphic a direct product
\[
J \times E,
\]
where $J$ is a a direct product of finitely many (possibly no) just infinite profinite groups and $E$ is a direct product of finite simple groups.  It follows that $J \times E$ is monomial; clearly, for this to be the case we must have either $J$ trivial or $E$ finite.  If $J$ is trivial, then $E$ is monomial, and then (LS) follows by the first paragraph.  So let us suppose that $E$ is finite.  Then by Lemma~\ref{lem:monomial:ji:bis} we know that $J$, and hence $G$, is locally isomorphic to a just infinite group $J'$.  If $J'$ is virtually abelian then (JA) follows.  If $J'$ is not virtually abelian, then $\QZ(J')=\triv$, from which it follows that $\QZ(G)$ is discrete.  There are then two possibilities as in Theorem~\ref{intro:wiltri}: if $\lnorm(G)$ is finite, then $J'$ is virtually a direct product of $n$ copies of a hereditarily just infinite profinite group, and we see that $2^n = |\lnorm(G)|$; after passing to an open subgroup, we obtain such a direct product that also occurs as a compact open subgroup of $G$, and hence (JH) holds.  If instead $\lnorm(G)$ is infinite, then $J'$ is a just infinite branch group, which means that $J'$ acts faithfully on $\lnorm(J')$ and hence the kernel of the action of $G$ on $\lnorm(G)$ is discrete.  At the same time, $\QZ(G)$ acts trivially on $\lnorm(G)$ and is the largest discrete normal subgroup of $G$; thus we have a faithful action of $G/\QZ(G)$ on $\lnorm(G)$, and (JB) follows.   Thus (i) implies (ii).

Conversely, suppose that (ii) holds.  In case (LS), it follows from the first paragraph that $G$ is locally monomial, and from Theorem~\ref{thm:propld} that $G$ has~\propLD.  For (JA), (JH) and (JB), we see that $G$ is locally isomorphic to a just infinite profinite group, and the same conclusions follow by Lemma~\ref{lem:monomial:ji:bis} and Theorem~\ref{thm:propld}.  Thus (ii) implies (i).
\end{proof}

Moving away from the locally monomial case, we now obtain a canonical factorization of an open subgroup.

\begin{defn}
Given a profinite group $H$ with~\propLD, a \defbold{monomial factorization} of $H$ is a factorization of $H$ as a direct product
\[
H = H_1 \times H_2 \times \dots \times H_n \times \prod_{S \in \mc{C}}H_S \times H_{\infty}
\]
with the following properties:
\begin{enumerate}[(i)]
\item $H_i$ is virtually just infinite for $1 \le i \le n$, and for $i \neq j$ then $H_i$ and $H_j$ are not similar;
\item $\mc{C}$ is a collection of isomorphism types of finite simple groups (possibly empty), and for each $S \in \mc{C}$ then $H_S$ is the direct product of $\aleph_0$ copies of $S$;
\item $H_{\infty}$ is either trivial or a direct product of infinitely many finite simple groups not belonging to $\mc{C}$, such that each isomorphism type only appears finitely many times.
\end{enumerate}
\end{defn}

Note in particular that the factors of the monomial factorization are monomial, with the exception of the leftover part $H_{\infty}$.

\begin{prop}\label{aji:reduced:new}
Let $G$ be a first-countable \tdlc group with~\propLD.
\begin{enumerate}[(i)]
\item $G$ has an open subgroup $H$ admitting a monomial factorization
\[
H = H_1 \times H_2 \times \dots \times H_n \times \prod_{S \in \mc{C}}H_S \times H_{\infty}.
\]
\item Let $N$ be an infinite compact monomial locally normal subgroup of $G$.  Then $N$ is virtually contained in one of the factors of the monomial factorization of $H$.
\item Let $K$ be a first-countable profinite group admitting a monomial factorization
\[
K = K_1 \times K_2 \times \dots \times K_m \times \prod_{S \in \mc{C}'}H_S \times K_{\infty},
\]
such that there is a local isomorphism $\phi: U \rightarrow H$ from $K$ to $H$.  Then the monomial factorizations are equivalent in the following sense:  We have $m = n$ and $\mc{C} = \mc{C}'$; there is a permutation $\pi$ of $\{1,\dots,m\}$ such that $[\phi(K_i \cap U)] = [H_{\pi(i)}]$; $[\phi(K_S)] = [H_S]$ for all $S \in \mc{C}$; and $[\phi(K_{\infty})] = [H_{\infty}]$.
\end{enumerate}
\end{prop}

\begin{proof}
(i)
By Theorem~\ref{thm:propld}, we can take
\[
G = \prod_{i \in I}L_i,
\]
where finitely many factors, say $L_1,\dots,L_{n'}$, are virtually just infinite profinite groups and the remaining factors are finite simple groups.  Since $G$ is first-countable, we can take the indexing set $I$ to be countable.  This immediately yields a factorization
\[
G = L_1 \times \dots \times L_{n'} \times \prod_{S \in \mc{C}}H_S \times H_{\infty},
\]
where $\mc{C}$ is the set of isomorphism types of finite simple factors of $G$ that occur infinitely many times, and $H_{\infty}$ is the product of the finite simple factors of $G$ whose isomorphism type only occurs finitely many times.  We can then group the factors $L_1,\dots,L_{n'}$ into similarity classes; by Lemma~\ref{lem:monomial:ji:bis}, the product of the factors in a given similarity class is virtually just infinite.  This yields the desired factorization of a group $H$ that has an open subgroup in common with $G$.

(ii)
By Theorem~\ref{thm:monomial}, after replacing $N$ with a finite index subgroup, we may assume either that $N$ is a direct product of $\aleph_0$ copies of a finite simple group $S$ or $N$ is virtually just infinite.

If $N$ is a direct product of $\aleph_0$ copies of a finite simple group $S$, let $\mc{I}_S$ be the set of elements of $\lnorm(G)$ with representatives locally isomorphic to $N$.  We see that $S \in \mc{C}$, $[N] \in \mc{I}_S$, and $[H_S]$ is the unique largest element of $\mc{I}_S$: we see the last asssertion, for example, by observing that $H/H_S$ has no locally normal subgroup locally isomorphic to $H_S$.  Thus $[N] \le [H_S]$, that is, $N$ is virtually contained in $H_S$.  If instead $N$ is virtually just infinite, one sees similarly that the set $\mc{J}$ of elements of $\lnorm(G)$ with representatives similar to $N$ has a unique largest element, which is represented by $H_i$ for $1 \le i \le n$, and hence $N$ is virtually contained in $H_i$.

(iii)
Writing $\theta([L]) = [\phi(L \cap U)]$, we have an order isomorphism $\theta$ from $\lnorm(K)$ to $\lnorm(H)$, which also preserves commensurability of the representatives.  In particular, applying (ii), we see that if $S \in \mc{C}'$, then $S \in \mc{C}$ and $\theta([K_S]) \le [H_S]$; after applying the same argument to $\theta\inv$, we see that $\mc{C} = \mc{C}'$ and $\theta([K_S]) = [H_S]$ for every $S \in \mc{C}'$.  Similarly, given $1 \le i \le m$ we see that $\theta([K_i]) \le [H_j]$ for some $1 \le j \le n$, ensuring that $H_j$ is similar to $K_i$, and conversely given $1 \le j \le n$ then $\theta\inv([H_j]) \le \theta([K_i])$ for some $1 \le i \le n$ such that $K_i$ is similar to $H_i$.  Thus in fact there is a permutation $\pi$ of $\{1,\dots,m\}$ such that $\theta([K_i]) = [H_{\pi(i)}]$, where $\pi$ is uniquely determined by the correspondence between the similarity classes of the virtually just infinite factors in the monomial factorizations of $K$ and $H$.  Finally, we can characterize $[K_{\infty}]$ in $\lnorm(K)$ as follows, and similarly for $[H_{\infty}]$ in $\lnorm(H)$: Given a closed locally normal subgroup $N$ of $K$ such that the finite simple normal subgroups of $N$ that do not belong to $\mc{S}$ generate a dense subgroup of $N$, then $N$ is virtually contained in $K_{\infty}$, and $[K_{\infty}]$ is the smallest element of $\lnorm(K)$ with this property.  Thus $\theta([K_{\infty}]) = [H_{\infty}]$.
\end{proof}

\begin{rem}
In Proposition~\ref{aji:reduced:new} we restricted to the first-countable case to avoid some technicalities, and because dropping this condition does not generalize the class of groups under consideration in any interesting way.  There are only countably many isomorphism types of finite group and every just infinite profinite group is first-countable, so the only difference between first-countable and general profinite groups with~\propLD in the context of Theorem~\ref{thm:propld} is that in the latter case, the factorization $\prod_{i \in I}L_i$ can include uncountably many copies of the same finite simple group.
\end{rem}

\subsection{Monomial constituents}\label{sec:monomial_constituents}

Given a \tdlc group $G$, a similarity class $\mc{S}$ of profinite groups is a \defbold{monomial constituent} of $G$ if groups in $\mc{S}$ are monomial and there exists $H \in \mc{S}$ such that $H$ is isomorphic to a locally normal subgroup of $G$.  Proposition~\ref{aji:reduced:new} shows that if $G$ has~\propLD, then any monomial factorization of a profinite group locally isomorphic to $G$ will yield all the monomial constituents of $G$.

Let us prove some analogues to Proposition~\ref{aji:reduced:new} describing normal subgroups of $G$ as a whole.  In particular, it will follow that the structure theory of \tdlc groups with~\propLD can to some extent be reduced to the locally monomial case.

\begin{prop}\label{prop:monomial_constituent}
Let $G$ be a first-countable \tdlc group with~\propLD and let $\mc{S}$ be a monomial constituent of $G$.  Then there is a continuous injective homomorphism $\phi: G_{\mc{S}} \rightarrow G$, where $G_{\mc{S}}$ is a first-countable \tdlc group, with the following properties:
\begin{enumerate}[(i)]
\item $G_{\mc{S}}$ is generated by its compact open subgroups, all of which belong to $\mc{S}$.
\item Given a compact locally normal subgroup $N$ of $G$ such that $N \in \mc{S}$, then $\phi$ restricts to a homeomorphism from $\phi\inv(N)$ to $N$.
\item $\phi(G_{\mc{S}})$ is characteristic in any open subgroup of $G$ that contains it.
\end{enumerate}
Moreover, given an open subgroup $K$ of $G$ containing $\phi(G_{\mc{S}})$, then the map $G_{\mc{S}} \rightarrow K$ has the same properties with respect to $K$ in place of $G$.
\end{prop}

\begin{proof}
Given a compact open subgroup $U$ of $G$, by Proposition~\ref{aji:reduced:new} there is an open subgroup of $U$ with a monomial factorization, such that exactly one of the factors $U_{\mc{S}}$ belongs to the similarity class $\mc{S}$.  

Now let $\mc{L}$ be the set of all compact locally normal subgroups of $G$ that belong to the similarity class $\mc{S}$.  Given $N \in \mc{L}$, we see by Proposition~\ref{aji:reduced:new}(ii) that $N$ is virtually contained in $U_{\mc{S}}$, and by Proposition~\ref{aji:reduced:new}(iii), $N$ commensurates $U_{\mc{S}}$.  It follows that the subgroup $\grp{\mc{L}}$ can be equipped with a \tdlc group topology such that $U_{\mc{S}}$ is embedded as an open subgroup; let $G_{\mc{S}}$ be $\grp{\mc{L}}$ with this topology and let $\phi$ be the inclusion map into $G$.  Note that $U_{\mc{S}} \in \mc{L}$, and up to finite index, $U_{\mc{S}}$ does not depend on the choice of $U$, so $G_{\mc{S}}$ does not depend on $U$.

Since $G$ is first-countable, so is $U_{\mc{S}}$, and hence $G_{\mc{S}}$ is first-countable.  It is clear from the construction that $\phi(G_{\mc{S}})$ is characteristic in $G$ and that (ii) holds; in particular, since $U_{\mc{S}}$ is open in $G_{\mc{S}}$, we see from (ii) that $\phi$ is continuous at the identity.  Since $\phi$ is evidently also a group homomorphism, it is continuous everywhere.  The profinite group $U_{\mc{S}}$, and hence all profinite groups commensurable with it, belong to $\mc{S}$; thus (i) holds.

Finally, if $K$ is an open subgroup of $G$ containing $\phi(G_{\mc{S}})$, then we see that the set of monomial factors in $\mc{S}$ of compact locally normal subgroups of $K$ belonging to $\mc{S}$ is the same as for $G$.  Thus we can identify $G_{\mc{S}}$ with $K_{\mc{S}}$ in the obvious way, ensuring that (iii) holds and the map $G_{\mc{S}} \rightarrow K$ has the same properties with respect to $K$ in place of $G$.
\end{proof}

We now prove the monomial decomposition theorem from the introduction.

\begin{proof}[Proof of Lemma~\ref{lem:monomial:qz}]
By Proposition~\ref{aji:reduced:new}, there is a compact open subgroup $U$ of $G$ with a factorization
\[
U = U_1 \times U_2 \times \dots \times U_n \times U_Q
\]
such that $U_Q = \ol{\QZ(U)}$; the factors $U_1,\dots,U_n$ are each virtually just infinite, with trivial quasi-centre; and no two factors are similar.  In particular, letting $Q_1 = \ol{\QZ(G)}$, we see that $Q_1 \cap U = U_Q$.  The quotient $G/Q_1$ then has a compact open subgroup of the form
\[
U^{*} = U_1 \times U_2 \times \dots \times U_n,
\]
which has trivial quasi-centre; thus $\QZ(G/Q_1)$ is discrete and $G/Q_1$ has~\propLD.  Taking $Q_2/Q_1 = \QZ(G/Q_1)$, we then obtain a characteristic quotient $G/Q_2$ of $G$ with~\propLD and trivial quasi-centre.  Notice that $\QZ(Q_2) \ge \QZ(G)$, so $\ol{\QZ(Q_2))}$ is open in $Q_2$.
\end{proof}

\begin{proof}[Proof of Theorem~\ref{thm:monomial:global}]
By Lemma~\ref{lcent:ldlat}, $G$ is [A]-semisimple; by Lemma~\ref{lem:melfin}, $G$ is first-countable.  We retain the compact open subgroup $U$ and its direct factors $U_1,\dots,U_n$ from the proof of Lemma~\ref{lem:monomial:qz}.  Let $\alpha_i$ be the element of $\lnorm(G)=\ldlat(G)$ represented by $U_i$ and let $M_i$ be the corresponding element $\CC^2_G(U_i)$ of $\LDlat(G)$; in particular, $M_i$ is closed and is the largest representative of $\alpha_i$.  Since $\{\alpha_1,\dots,\alpha_n\}$ forms a partition of $\infty$ in $\ldlat(G)$, we have an open subgroup
\[
M := \grp{M_1,M_2,\dots,M_n} = M_1 \times M_2 \times \dots \times M_n.
\]
Let $N$ be a closed locally normal subgroup of $G$ that is locally similar to $J_i$.  By Proposition~\ref{aji:reduced:new}(ii), $U_i \cap N$ is open in $N$, so also $M_i \cap N$ is open in $N$.  By Theorem~\ref{thm:cent_pseudocomp} it follows that $M_j$ centralizes $N$ for all $j \neq i$, and hence $N \le M_i$.  Thus (ii) holds.

Consider a closed locally normal subgroup $K$ of $G$ and let $\beta_i$ be the locally normal subgroup of $K$ represented by $K \cap U_i$; note that the product $\prod^n_{i=1}(K \cap U_i)$ is a compact open subgroup of $K$, and also that $K$ has~\propLD and trivial quasi-centre by Lemmas~\ref{lem:ld:hered} and~\ref{lem:ld:qz}.  For $1 \le i \le n$, since the compact representatives of $\alpha_i$ are monomial, they are similar to the representatives of $\beta_i$.  Since the representatives of $\alpha_i$ lie in different similarity classes for $1 \le i \le n$, we deduce that $\prod^n_{i=1}(K \cap U_i)$ is a monomial factorization of an open subgroup of $K$, and by Proposition~\ref{aji:reduced:new}(iii), it follows that every automorphism of $K$ will preserve each of the elements $\beta_1,\dots,\beta_n$ of $\ldlat(K)$.  In particular, if $M_i \le K$ for some $i$, we see that $\alpha_i = \beta_i$ and that $M_i$ is the largest representative of $\beta_i$, so $M_i$ is characteristic in $K$.  This proves (i).

Suppose $G$ contains a compact just infinite representative $J_i$ of $[M_i]$.  Then $J_i \le M_i$ by (ii), and hence $\lnorm(M_i)^{M_i} = \{0,\infty\}$, ensuring that every closed normal subgroup of $M_i$ is discrete or open in $M_i$.  The former is ruled out by the fact that $\QZ(M_i)=\triv$; thus (iii) holds.

Finally, let us suppose that $G = \ms{L}(G)$.  Given $g \in G$, then conjugation by $g$ induces a local automorphism of $U$; using Proposition~\ref{aji:reduced:new}(iii), we can restrict this local automorphism to an isomorphism
\[
V_1 \times V_2 \times \dots \times V_n \rightarrow gV_1g\inv \times gV_2g\inv \times \dots \times gV_ng\inv,
\]
where $V_i$ and $gV_ig\inv$ are open subgroups of $U_i$.  We then see that there are elements $g_1,\dots,g_n$ of $G$ such that $g_ivg\inv_i =gvg\inv$ for all $v \in V_i$, but $g_i \le \CC_G(V_j)$ for all $j \neq i$.  We then have $g_i \in M_i$ and $g = g_1g_2 \dots g_n$.  This proves that $G = M$.  On the other hand, we see that every local automorphism of $M_i$ extends to a local automorphism of $G$, and hence can be realized by conjugation in $G$; thus $M_i = \ms{L}(M_i)$ for $1 \le i \le n$.  Finally, it is now clear that $J_i$ appears as a compact open subgroup of $M_i$, completing the proof of (iv).  
\end{proof}

\begin{rem}
In the situation of Theorem~\ref{thm:monomial:global}, one sees that factors $M_1,\dots,M_n$ correspond to the monomial constituents $\mc{S}_1,\dots,\mc{S}_n$ of $G$ from Proposition~\ref{prop:monomial_constituent}, and for $1 \le i \le n$ the map $\phi: G_{\mc{S}_i} \rightarrow G$ constructed in Proposition~\ref{prop:monomial_constituent} restricts to an open embedding of $G_{\mc{S}_i}$ into $M_i$.
\end{rem}

\subsection{Groups of type (JH)}\label{sec:hji}

We now give the analogue of Theorem~\ref{intro:bew} for groups of type (JH).

\begin{thm}\label{thm:hji:atomic}
Let $G$ be a \tdlc group of type (JH), with $\QZ(G)=\triv$.  Then there is a hereditarily just infinite profinite group $J$ and an open subgroup
\[
P = P_1 \times P_2 \times \dots \times P_n
\]
of $G$, where the factors $P_1,\dots,P_n$ (which we call the \defbold{atomic factors} of $G$) have the following properties:
\begin{enumerate}[(i)]
\item Every automorphism of $G$ permutes the factors $\{P_1,\dots,P_n\}$ (possibly trivially); in particular, $P$ is characteristic in $G$.
\item For $1 \le i \le n$, $\QZ(P_i) = \triv$ and $P_i$ has a compact open subgroup isomorphic to $J$.
\item Given $1 \le i \le n$, exactly one of the following holds:
\begin{description}
\item[(Reducible type)] $P_i \cap \Res(G) = \triv$;
\item[(Mysterious type)] $P_i \cap \Res(G)$ is open, but $P_i \cap \Res(\Res(G)) = \triv$;
\item[(Simple type)] $P_i \cap \Res(G) = \Res(P_i)$ is open in $P_i$ and topologically simple.
\end{description}
\end{enumerate}
\end{thm}

\begin{proof}
By Theorem~\ref{thm:monomial}, $\lnorm(G)$ is finite and there is an open subgroup $J'$ of $G$ that is a direct product
\[
J' = J_1 \times J_2 \times \dots \times J_n,
\]
where $J_i \cong J_j$ for all $i$ and $j$ and $J:= J_1$ is hereditarily just infinite.  Note that the elements $[J_1],\dots,[J_n]$ of $\lnorm(G)$ are distinct and are exactly the atoms of $\lnorm(G)$.  We set $P_i = \QC^2_G([J_i])$; then we have an open subgroup
\[
P = P_1 \times P_2 \times \dots \times P_n.
\]
Since the factors $P_1,\dots,P_n$ are obtained from $\lnorm(G)$ in a canonical way, they form a characteristic class of subgroups; thus (i) is satisfied.  We also note that $\QZ(P_i) \le \QZ(G) = \triv$ for all $1 \le i \le n$, and clearly $J_i \le P_i$; this proves (ii).

It remains to divide the atomic factors into three types as in (iii).  Without loss of generality, we can replace $G$ with the finite index open subgroup $\bigcap^n_{i=1}\N_G(P_i)$; note that this does not change the fact that $G$ is compactly generated, nor does it change $\Res(G)$.  Clearly the three listed types are mutually exclusive; note also that because of (ii), any closed locally normal subgroup of $P_i$ is either trivial or open.  Fix $1 \le i \le n$ and suppose $P_i$ is not of reducible or mysterious type.  Then we see that $\Res(\Res(G)) \cap P_i$ must be open in $P_i$; since $P$ is open and characteristic in $G$, we see that $\Res(G) \le P$ and hence $P_i \cap \Res(\Res(G)) = \Res(R)$, where $R = P_i \cap \Res(G)$.  In particular, the product of $\Res(R)$ with $P_j$ for $j \neq i$ is an open normal subgroup $O$ of $G$; we then have
\[
\Res(R) = P_i \cap O \ge P_i \cap \Res(G) = R,
\]
thus $\Res(R) = R$.  Since in addition, every nontrivial closed normal subgroup of $R$ is open, we conclude that $R$ is topologically simple.  We then have 
\[
\Res(R) \le \Res(P_i) \le R,
\]
so $R = \Res(P_i)$.  Thus $P_i$ is of simple type, completing the proof.
\end{proof}

We now prove the theorem from the introduction about mysterious type.  Part (iv) is derived from a more general fact, which we prove separately in a lemma.

\begin{lem}[{\cite[Corollary~1.4]{ReidE}}]\label{endolem}Let $U$ be a profinite group such that $U$ has only finitely many open subgroups of each index, and such that $U$ is isomorphic to a proper open subgroup of itself.  Then $U$ has an infinite abelian normal subgroup.\end{lem}

\begin{lem}\label{uniext}Let $G$ be a \tdlc group.  Suppose that $G = \langle K,x\rangle$, where $K$ is a residually discrete open normal subgroup of $G$, and suppose that $G$ has no nontrivial fixed points in its action on $\lnorm(G)$.  Suppose also that $K$ has no infinite abelian locally normal subgroups and that every compact open subgroup of $K$ has only finitely many open subgroups of each index.  Then $G$ is residually discrete; indeed, every open normal subgroup of $K$ contains an open normal subgroup of $G$.\end{lem}

\begin{proof}Let $H$ be an open normal subgroup of $K$ and let $U$ be a compact open subgroup of $H$.  Define the following subgroups:
\[ H_+ = \bigcap^\infty_{i=0} x^iHx^{-i}; \; H_- = \bigcap^\infty_{i=0} x^{-i}Hx^i; \; U_+ = \bigcap^\infty_{i=0} x^iUx^{-i}; \; U_- = \bigcap^\infty_{i=0} x^{-i}Ux^i.\]
Note that all conjugates of $H$ and their intersections are closed and normal in $K$; furthermore, $xH_+x^{-1}$ and $x^{-1}H_+x$ are both locally equivalent to $H_+$.  Hence $\alpha_+ = [H_+]$ is an element of $\lnorm(G)$ that is fixed by $G$.  It follows that $\alpha_+ \in \{0,\infty\}$.  By \cite[Lemma~1]{Willis}, the set $U_+U_-$ is a neighbourhood of the identity in $H$.  Suppose that $\alpha_+=0$; then $H_+$ is discrete, so $U_+$ is finite, and hence $U_-$ is open in $K$.  By construction, $xU_-x^{-1}$ is a proper open subgroup of $U_-$ that is isomorphic to $U_-$.  Thus $K$ has an infinite abelian locally normal subgroup by Lemma~\ref{endolem}, a contradiction.  Thus $\alpha_+ = \infty$, in other words, $H_+$ is open in $K$.  By the same argument, $H_-$ is also open in $K$, so $H_+ \cap H_-$ is open in $K$.  Hence every open normal subgroup of $K$ contains an open normal subgroup of $G$.
\end{proof}

\begin{proof}[Proof of Theorem~\ref{thm:hji:mysterious}]
Let $U$ be a proper compact open subgroup of $R$ and note that the $G$-conjugates of $U$ have trivial intersection.  By \cite[Theorem~1.2]{ReidEqui}, there is some $x \in R \setminus \{1\}$ such that the $G$-conjugacy class of $x$ accumulates at the identity.  In particular, every open subgroup of $R$ contains $gxg\inv$ for some $g \in G$.

Now let $K = \cgrp{rxr\inv \mid r \in R}$ and consider a nontrivial closed normal subgroup $Q$ of $R$.  Then $K$ and $Q$ are both open, because $R$ is locally \hji and has trivial quasi-centre.  We then see that $Q$ contains $gxg\inv$ for some $g \in G$; since $Q$ is normal in $R$, it follows that $Q$ contains $rgxg\inv r\inv$ for every $r \in R$.  Since $R$ is normal in $G$, in fact $Q$ contains $grxr\inv g\inv$ for every $r \in R$; since $Q$ is a closed subgroup of $R$, it follows that $gKg\inv \le Q$.  We have now proved (i) except for the fact that $K$ is not compact.

Consider the case $Q = K \cap gKg\inv$, where $g \in G$ is such that $K \nleq gKg\inv$.  Then $Q$ contains $hKh\inv$ for some $h \in G$, and then we see that $hKh\inv < K$.  This proves (ii).  However, by Lemma~\ref{endolem}, no compact open subgroup of $R$ is isomorphic to a proper open subgroup of itself.  Thus $K$ is not compact, completing the proof of (i).

We note by (i) that $R$ has no compact open normal subgroups.  It follows that if $H$ is a compactly generated subgroup containing $R$, then $H$ also has no compact open normal subgroups.  By \cite[Corollary~4.1]{CM} it follows that $\Res(H) \neq \triv$, and hence $\Res(H)$ is open.  Since $R$ is open and normal in $H$, we have $\Res(H) = \Res_R(H)$.  Since $\Res(R) = \triv$, we see that $R$ cannot be compactly generated, proving (iii).

Finally, consider a subgroup $L$ of $G$ containing $R$ such that $L/R$ is virtually polycyclic; then there is $R \le L_0 \le L$ such that $L_0/R$ is polycyclic and $L_0$ is a normal subgroup of $L$ of finite index.  By repeated application of Lemma~\ref{uniext} we see that every open normal subgroup of $R$ contains an open normal subgroup of $L_0$.  In turn, every open normal subgroup $M$ of $L_0$ contains the open normal subgroup $M' = \bigcap_{y \in Y}yMy\inv$ of $L$, where $Y$ is a finite set of coset representatives for $L_0$ in $L$.  This proves (iv).
\end{proof}

\subsection{Groups locally isomorphic to just infinite profinite branch groups}

We now turn to groups $G$ of type (JB).  In this case there is no atomic decomposition coming from the structure lattice, since $\lnorm(G)$ is atomless.  Instead, we take the approach of decomposing $\Res(G)$ into directly indecomposable parts.  Other than that, we obtain a statement similar to Theorem~\ref{thm:hji:atomic}.

\begin{thm}\label{thm:locaji:structure}
Let $G$ be a \tdlc group of type (JB), with $\QZ(G)=\triv$.  Then $\Res(G)$ is a direct factor of an open subgroup and admits a decomposition into finitely many directly indecomposable direct factors
\[
\Res(G) = P_1 \times P_2 \times \dots \times P_n,
\]
which we call the \defbold{components} of $G$, with the following properties.
\begin{enumerate}[(i)]
\item Every direct factor of $\Res(G)$ is a direct product of a subset of the components.
\item Each component $P$ is noncompact, but locally isomorphic to a just infinite profinite branch group, and we have $\lnorm(P)^P = \{0,\infty\}$.  Moreover, $G$ does not normalize any proper nontrivial closed subgroup of $P$.
\item Given $1 \le i \le n$, exactly one of the following holds:
\begin{description}
\item[(Mysterious type)] $P_i$ is residually discrete;
\item[(Simple type)] $P_i$ is topologically simple.
\end{description}
\end{enumerate}
\end{thm}

\begin{proof}
We see that $R = \Res(G)$ is itself of type (JB) with trivial quasi-centre.  In particular, we have $\lnorm(R) = \ldlat(R)$, $R$ is [A]-semisimple, and $R$ acts faithfully on $\lnorm(R)$ with finitely many fixed points.  Let $\alpha_1,\dots,\alpha_n$ be the atoms of $\lnorm(R)^R$ and set $P_i = \QC^2_R(\alpha_i)$ for $1 \le i \le n$.  Then we obtain an open subgroup of $R$
\[
O =  P_1 \times P_2 \times \dots \times P_n.
\]
From the construction we see that $O$ is characteristic in $R$, hence normal in $G$; from there we see that $\CC_G(R) \times O$ is an open normal subgroup of $G$.  Thus $R \le \CC_G(R) \times O$ and in fact $R = O$.  It is now clear that $\lnorm(P_i)^{P_i} = \{0,\infty\}$ for $1 \le i \le n$; in particular, $P_i$ is directly indecomposable, and every nontrivial closed normal subgroup of $P_i$ is open.  Moreover, $P_i$ is of type (JB), so it is locally isomorphic to a just infinite profinite branch group.  Given a nontrivial closed subgroup $Q$ of $P_i$ that is normalized by $G$, we see that $Q$ is open in $P_i$, and hence the product $Q'$ of $Q$ with the components other than $P_i$ is an open normal subgroup of $G$.  But then $R \le Q'$, so $P_i \le Q$.  To complete the proof of (ii), suppose for a contradiction $P_i$ is compact.  Then $P_i$ is a just infinite profinite group; in particular, $\triv < \M(P_i) < P_i$.  But then $\M(P_i)$ is characteristic, hence $G$-invariant, and we have a contradiction to (ii).  Thus $P_i$ is noncompact as claimed.

Consider now a direct factor $D$ of $\Res(G)$.  Since $R$ is [A]-semisimple, $D$ has a unique direct complement $\CC_R(D)$, and then $D = \CC^2_R(D)$, so $D$ is closed and is the largest representative of an element $\alpha \in \lnorm(R)$.  Since $R = D \times \CC_R(D)$, in fact $\alpha \in \lnorm(R)^R$, so we can write $\alpha$ as the join of $\{\alpha_i \mid i \in I\}$ for some subset $I$ of $\{1,\dots,n\}$.  Given the direct decomposition of $R$ into components, we see that in fact $D$ is generated by the components $\{P_i \mid i \in I\}$.  This proves (i).  

It remains to divide the components into two types.  Fix $1 \le i \le n$.  If $\Res(P_i) > \triv$, then by (ii), we must have $\Res(P_i) = P_i$.  Since every nontrivial closed normal subgroup of $P_i$ is open, it follows that $P_i$ is topologically simple.  Otherwise, clearly $P_i$ is of mysterious type.\end{proof}

Similar to type (JH), there are no known examples of mysterious components; however, unlike for type (JH), we can rule out the mysterious components when $G$ is compactly generated.  In fact, in the present situation we can prove a result about $\Res_G(H)$, where $H$ is any compactly generated subgroup of $G$.  We first recall a sufficient condition for a \tdlc group to have a nontrivial contraction group.

\begin{lem}[{See \cite[Proposition~6.14 and Theorem~6.19]{CRW-Part2}}]\label{lem:lcent_contraction}
Let $G$ be a compactly generated \tdlc group that is [A]-semisimple and let $\mc{A}$ be a subalgebra of $\lcent(G)$ on which $G$ acts faithfully.
\begin{enumerate}[(i)]
\item Suppose that $G$ has a compact open subgroup $U$ such that $\bigcap_{g \in G}gUg\inv = \triv$.  Then there is a finite subset $\{\alpha_1,\dots,\alpha_n\}$ of $\mc{A}$ such that for all $\beta \in \mc{A} \setminus \{0\}$, there is some $g \in G$ and $i \in \{1,\dots,n\}$ such that $g\alpha_i < \beta$.
\item Let $V$ be a compact open subgroup of $G$, and suppose there is $g \in G$ and $\alpha \in \mc{A}$ such that $g\alpha < \alpha$. Then there is a natural number $n_0$ such that
\[
\QC^2_U(g^{n_0}\alpha \setminus g^{n_0+1}\alpha) \le \con(g).
\]
\end{enumerate}
\end{lem}

\begin{lem}\label{lem:branch:simpleres}
Let $G$ be a \tdlc group of type (JB), with $\QZ(G)=\triv$, and let $H$ be a compactly generated subgroup of $G$.  Then $\Res_G(H)$ is the direct product of finitely many (possibly none) topologically simple groups, each of which is locally normal in $G$.  Moreover, if $\Res_G(H) \le H$, then
\[
\Res_G(H) = \ol{G^\dagger_H} = \ol{H^\dagger}.
\]
\end{lem}

\begin{proof}
By Lemma~\ref{lem:reduced_envelope} and replacing $G$ with an open subgroup, we may assume that $\Res_G(H) = \Res(G)$.  Since $\Res(G)$ is a direct factor of an open subgroup, there is a natural isomorphism of topological groups 
\[
\Res(G) \rightarrow \Res(G)\CC_G(\Res(G))/\CC_G(\Res(G)) = \Res(G/\CC_G(\Res(G))).
\]
If $\Res(G)=\triv$ there is nothing more to prove, so without loss of generality, $\Res(G)$ is locally similar to $G$.  We can now replace $G$ with $G/\CC_G(\Res(G))$ and hence assume $\CC_G(\Res(G)) = \triv$.  As a result, $\Res(G)$ is open in $G$.

By Theorem~\ref{thm:locaji:structure}, we have
\[
\Res(G) = P_1 \times P_2 \times \dots \times P_n,
\]
where $P_1,\dots,P_n$ are the components of $G$; moreover, we see that the elements $\alpha_i = [P_i]$ of $\ldlat(G)$ are exactly the atoms of $\ldlat(G)^G$, while $P_1,\dots,P_n$ are exactly the minimal nontrivial closed normal subgroups of $G$.  Applying Lemma~\ref{lem:lcent_contraction}, for $1 \le i \le n$ there is $g_i \in G$ and $0 < \beta_i < \alpha_i$ such that $g_i\beta_i < \beta_i$, and then the contraction group $\con_{P_i}(g_i)$ of $g_i$ acting on $P_i$ is nontrivial.  In particular, the semidirect product $P_i \rtimes \grp{g_i}$ is not residually discrete.  Since $P_i$ acts on $\ldlat(P_i) = \lnorm(P_i)$ with no nontrivial fixed points, it follows from Lemma~\ref{uniext} that $P_i$ is not residually discrete.  Thus $P_i$ cannot be of mysterious type, so by Theorem~\ref{thm:locaji:structure} it is topologically simple.

In the case that $\Res_G(H) \le H$, it is clear that $\ol{G^\dagger_H} = \ol{H^\dagger}$, and the existence of the elements $g_i$ makes it clear that $\ol{H^\dagger} = \Res_G(H)$.
\end{proof}

\subsection{Noncompact topologically simple local direct factors}\label{sec:simple}

We can now prove the last two theorems from the introduction.

\begin{proof}[Proof of Theorem~\ref{thm:simpleres}]
By Lemma~\ref{lem:reduced_envelope} and replacing $G$ with an open subgroup, we may assume that $\Res_G(H) = \Res(G)$.  We first take the monomial factorization of an open normal subgroup of $G$,
\[
M = M_1 \times M_2 \times \dots \times M_m,
\]
as in Theorem~\ref{thm:monomial:global}; it is then clear that $\Res_G(H)$ is the direct product of $\Res_{M_i}(H)$ for $1 \le i \le m$.  Thus we reduce to the case when $G$ is locally monomial, and hence of type (JH) or (JB).

If $G$ is of type (JH), we are in the case that $H\Res_G(H)/\Res_G(H)$ is virtually polycyclic.  We take the atomic factorization of $G$ as in Theorem~\ref{thm:hji:atomic}, and see by Theorem~\ref{thm:hji:mysterious}(iv) that only the topologically simple atomic factors of $G$ can contribute to $\Res_G(H)$, and hence there are no atomic factors of mysterious type; in other words, $\Res_G(H)$ is a finite direct product of topologically simple groups.

If $G$ is of type (JB), then the conclusions follow from Lemma~\ref{lem:branch:simpleres}.
\end{proof}

\begin{proof}[Proof of Theorem~\ref{thm:aji:titscore}]
We see that $\ol{G^\dagger}$ is a closed $G$-invariant subgroup of $\Res(G)$, and by Theorem~\ref{thm:monomial:global}, $\Res(G)$ splits into monomial parts.  Given Theorems~\ref{thm:hji:atomic} and~\ref{thm:locaji:structure}, we can write
\[
\Res(G) = P_1 \times P_2 \times \dots \times P_n,
\]
where $\{P_1,\dots,P_n\}$ is a characteristic class of subgroups of $G$, and each of the factors $P_i$ is one of the following:
\begin{enumerate}[(a)]
\item the intersection of an atomic factor of mysterious type with $\Res(G)$;
\item a locally branch component of mysterious type;
\item a topologically simple group that is locally isomorphic to a just infinite profinite group.
\end{enumerate}
Note that in all cases, $\lnorm(P_i)^{P_i} = \{0,\infty\}$.

After passing to a finite index open subgroup (which changes neither $G^\dagger$ nor $\Res(G)$), we may assume $P_i$ is normal in $G$ for $1 \le i \le n$.  We then see that
\[
\ol{G^\dagger} = Q_1 \times Q_2 \times \dots \times Q_n,
\]
where $Q_i = \ol{(P_i)^\dagger_G}$.

Fix $1 \le i \le n$ and suppose $Q_i \neq \triv$; that is, there is some $g \in G$ with a nontrivial contraction group on $P_i$.  In particular, the group $P_i \rtimes \grp{g}$ is not residually discrete.  By Lemma~\ref{uniext}, it follows that $P_i$ is not residually discrete, so we are in case (c), that is, $P_i$ is topologically simple.  In particular, $P_i = Q_i$.

Now suppose $G$ is compactly generated.  The group $\Res(G)/\ol{G^\dagger}$ is isomorphic to the direct product of those $P_i$ such that $Q_i =\triv$.  By Lemma~\ref{lem:branch:simpleres}, the factors $P_i$ of type (b) are ruled out, and for type (c), if $P_i$ is locally branch we have $P_i = Q_i$.  Thus $\Res(G)/\ol{G^\dagger}$ is a direct product of finitely many groups of type (JH); in particular, it is locally isomorphic to a finite direct product of \hji profinite groups.
\end{proof}

Theorem~\ref{thm:aji:titscore} has the following corollary, which illustrates a significant connection between topologically simple \tdlc groups and the internal structure of just infinite profinite groups.

\begin{cor}\label{cor:jiuni}Let $U$ be a just infinite profinite group that is not virtually abelian.  Let $V \le_o U$ and let $\theta: V \rightarrow U$ be a continuous injective open homomorphism such that $\theta(W) \not= W$ for all $W \le_o V$.  Then there is a compactly generated \tdlc group $G$ with an open normal topologically simple subgroup $S$, such that $S$ is locally similar to $H$.  Consequently, the composition factors of $U$ are of bounded order.
\end{cor}

\begin{proof}As explained in \cite{BEW}, $U$ naturally embeds in $L = \ms{L}(U)$ as a compact open subgroup, and then $\theta$ is the restriction of an inner automorphism of $L$, induced by the element $x$ say of $L$.  By the construction of the group of germs, $L$ is locally just infinite and $\QZ(L)=\triv$.  By the assumptions on $\theta$ and Lemma~\ref{lem:anisotropic}, we see that $L^\dagger_x \neq \triv$.  It follows that $H^\dagger \neq \triv$, where $H = \langle U, x \rangle$.  By Theorem~\ref{thm:aji:titscore}, $H$ has a topologically simple locally normal subgroup $S$.  Since $U$ is monomial, $S$ is locally similar to $U$.  Finally, we form the group $G = \N_H(S)/\CC_H(S)$, which has $S$ as an open normal subgroup.

By \cite[Proposition~4.6]{CRW-Part2}, taking a compact open subgroup $V$ of $S$, then the composition factors of $V$ are of bounded order.  Since $U$ is similar to $V$, it follows that the composition factors of $U$ are also of bounded order.
\end{proof}

\end{document}